\documentclass[10pt,a4paper]{amsart}
\usepackage{pictex}
\input amssym.def
\input amssym.tex
\title[Hecke groups ]
{arithmetic and geometry of the Hecke groups
}

\author{\SMALL Cheng Lien Lang}
\author{ \SMALL  mong lung lang}

\begin{document}

\baselineskip=12pt


\maketitle

\vspace{-0.3in}

\begin{abstract}
We study the arithmetic and geometry properties of the Hecke  group $G_q$.
In particular, we prove that  $G_q$ has a subgroup $X $ of index $d$, genus $g$ with
  $v_{\infty} $ cusps, and  $\tau_2$ (resp. $v_{r_i}$) conjugacy classes of
 elements that are conjugates of $S$ (resp. $R^{q/r_i}$)
 if and only if
(i)
$ 2g-2 + \tau_2/2 +\sum_{i=1}^k  v_{r_i}(1-1/r_i)  + v_{\infty} = d(1/2-1/q)$,
 and
(ii) $ m _0= 4g-4 +\tau_2 + 2 v_{\infty} + \sum _{i=1}^k v_{r_i}(2-q/r_i)\ge 0$ is a multiple of $q-2$, (iii) $m \ge 0$.
In the case $q$ is odd, (ii) is a consequence of (i).

\end{abstract}

\section{Introduction}
\subsection{}
The (inhomogeneous) Hecke group
 $G_q$ is defined to be the maximal discrete subgroup
of  $PSL (2,\Bbb R)$
 generated by  $S$ and $T$, where $\lambda _q =2$cos$\,(\pi/q)$,
 $$ S =
\left (
\begin{array}{rr}
0 & 1 \\
-1 & 0 \\
\end{array}
\right ) \,,\,\,
  T = \left (
\begin{array}{rc}
1 & \lambda  _ q\\
0 & 1 \\
\end{array}
\right ) \,.\eqno(1.1)
$$
Let $R= ST^{-1}$.
Then $R$ has order $q$
and $\{S, R\}$ is a set of independent generators of $G_q$.
 Equivalently, $G_q$ is a free product of $\left <S\right >$
  and $\left <R\right >$.
The main purpose of this article is to study the geometric and
 arithmetic properties of subgroups of finite index of $G_q$.

\subsection{}
 The set of cusps of $G_q$ is $\Bbb Q[\lambda_q] \cup \{\infty\}$ if and only
  if $q= 3, 5$. We will give an inductive procedure
  (induction on the depth of $q$-gons) that enables us to generate the set of
   cusps of $G_q$  (Lemma 3.2).
   As the index of $G_q$ in $PSL(2, \Bbb Z[\lambda_q])$ is infinite, it is
    important to characterise  members of $G_q$.
    A simple algorithm that determines whether a matrix
    of $PSL(2, \Bbb Z[\lambda_q])$ belongs to $G_q$ can be found in Proposition 3.7.
     The  algorithm  can be implemented in a computer.

\subsection{}
A set of generators $\{x_i\}$ of $X$ is called a set of {\em independent generators}
 if $X$ is a free product of the cyclic groups $\left <x_i\right >$.
  $G_q$ is a free product of $\left <S\right >$
  and $\left <R\right >$.
  By Kurosh's Theorem, every subgroup $X$ of finite index of $G_q$ has a set of independent generators.
 Proposition 4.4 and Theorem 5.2 demonstrate how  arithmetic and geometric properties can be combined to give an inductive procedure for finding a special polygon  (fundamental
  domain) $M_X$ and
 a set of independent generators $I_X$ for $X$
 (the case $q$ is a prime has been done in [LLT1]). In particular, this is applied
 to the principal congruence subgroup of level 2, the commutator subgroup $G_q'$
 and subgroups of index 2
  (subsection 5.4).

\subsection{} As a special case of the Hurwitz-Nielsen realisation problem,
 Millington [Mi] showed that as long as
 $d = 3\tau_2 +4v_3 +12q +6t-12$, then the modular group $G_3$ possesses
  a subgroup $X$ of index $d$, such that $X\setminus \Bbb H$ has $\tau_2$ (resp. $v_3$)
   elliptic points of order 2 (resp. 3),
    $t$ cusps, and genus $g$. We are able to generalise this result to $G_q$
    by studying the Hecke-Farey symbols   (see Section 6). As mentioned in [K1],
         the problem of recognising $X $ as a normal
         subgroup of certain geometric invariants has been  left aside in the literature. Our study of this topic starts
          with some elementary observation of the permutation representations
           of $S$ and $R$ on the set of cosets $G_q/X$ which we will elaborate
            more in subsection 1.5.

     \subsection{}
      We prove that
the  action of $S$ and $R$ on $G_q /X$ is isomorphic to their action on $M_X$
 (see Lemma 7.1) and that
 the permutation representations of $S$ and $R$ on $M_X$ can be obtained
  by a simple reading of the special polygon $M_X$
   ((7.4) and (7.9)).
    These two representations $f(S)$ and $f(R)$ carry some important information
    about $X$. In particular,
    \begin{enumerate}
    \item[(i)]
  $f(S)$ and $f(R)$ can be used to determine
     whether  $X$ is normal in $G_q$ and its normaliser in $PSL(2, \Bbb R)$ (Proposition 8.1 and Discussion 8.2),
    \item[(ii)] in the case $q=3$, $f(S)$ and $f(R)$ can be used to determine
     whether $X$ is congruence (see Section 11 and [H]),
    \item[(iii)]  $f(S)$ and $f(R)$ can be used to study
     the geometric invariants of $X$ (Section 9) and
      {\em Dessins d'enfants} (see pp.12 of [HR]). Propositions 9.2-9.4
       study the possible realisation of a group $Y$ as a normal subgroup
        of $G_q$.
     \end{enumerate}

\subsection{} It is well known that
there is a correspondence between the set of {\em maps}
 and the set of subgroups of finite index of $G_q$ and  that
 the maps are uniquely  determined
   by their  {\em map subgroups}
  (see [JS], [CS], [IS]).
  Let $X$ be a subgroup of $G_q$ of finite index.
  We give a detailed construction of the map $M(X)$ whose map subgroups
   are conjugates of $X$. Both $M(X)$ and $X$ are explicitly given.
      Aut$\,M(X)$   can be determined as well (see Section 10).

\section {Tessellation of the  upper half plane }
Let $D^*$ denote the $(2,q,\infty)$ triangle with vertices $i$, $e^{\pi i/q}$
 and $\infty$. $D^*$ is a fundamental domain of the Coxeter group
  $G_q^*$ generated by reflections along the sides of $D^*$.
   Hecke group $G_q$ is the subgroup of index 2 consists of
    all the orientation preserving isometries.

Let $\Bbb H$ be the union of the upper half plane and the set $\{ g(\infty)\,:\, g\in G_q^*\}$.
The $G_q^*$ translates of $D^*$ form a tessellation $\mathcal I^*$ of $\Bbb H$
 (endowed with the hyperbolic metric) by $(2,q, \infty)$ triangles.
  The $G_q^*$ translates of $i$, $e^{\pi i/q}$ and $\infty$
   are called {\em even vertices, odd vertices} and {\em cusps} ({\em free
    vertices}) of $\mathcal I^*$ respectively.
     The $G_q^*$ translates of the hyperbolic line joining $i$ to $\infty$ (resp. $e^{\pi i/q}$ to
     $\infty$) are called {\em even edges} (resp. {\em odd edges}) of $\mathcal I^*$.
        The $G_q^*$ translates of the hyperbolic line joining $i$ to $e^{\pi i/q}$
         are called {\em f-edges } of $\mathcal I^*$.
       The hyperbolic line $(0, \infty)$ consists of two even edges. The $G_q^*$ translates
        of $(0, \infty)$ are called the {\em even lines} of $\mathcal I^*$.
        {\em  The hyperbolic line joining $x$ and $y$ is denoted by $(x, y)$.}

   The set of even lines give a tessellation of $\Bbb H$ into
 {\em ideal q-gons}, that is, hyperbolic $q$-gons with $q$ cusps.
  Note that their vertex angle is 0. Each $q$-gon contains a unique odd vertex.

  The  $f$-edges form a $q$-regular tree, where the odd vertices are
  considered as the set of vertices of this tree. We introduce a vertex of
   valence two to this $q$-regular tree at $i$, denoted by $v_0$. Each $q$-gon $P$
    contains a unique vertex $v_P$ of this tree.
     The {\em  depth}  of $P$  denoted by
     $d(P)$ is defined to the
      the distance between $v_0$ and  $v_P$ (the distance between adjacent
      vertices is 1).

\section{Cusps and reduced forms }
 The main purpose of this section is to  determine  whether
  a matrix of $PSL(2, \Bbb Z[\lambda_q])$ belongs to $G_q$ (Proposition 3.7).
   Lemma 3.2 gives the set of  cusps of $G_q$.

\subsection{ Reduced forms of cusps of $G_q$}  The set of cusps of $G_q$ is a subset of $\Bbb Q(\lambda_q)\cup \{\infty\}$.
Let $x$ be a cusp. We say $x= a/b$ is in {\em reduced form} if
\begin{enumerate}
\item[(i)] there exists $c/d$ such that
{\tiny $\left (
\begin{array}{rr}
a & c \\
b & d \\
\end{array}
\right )$} $\in G_q$,  and (ii) $b\ge 0$.
\end{enumerate}

\smallskip
Let $g$ be given as in (i) of the above.
Since $w=$  {\tiny $\left (
\begin{array}{cr}
1 & 0 \\
0 & -1 \\
\end{array}
\right )$} normalises $G_q$, $ -wgw^{-1} \in G_q$. Hence $a/b$ is in reduced
 form if and only if
 $-a/b$ is in reduced form.   In the case $a, b \ge 0$, it is also easy to see
  that $a/b$ is in reduced form if and only if $-b/a$ is in reduced form (study $S^{-1}g$)
    if and only if $b/a$ is in reduced form
    (study $w(S^{}g)w^{-1}$).

\subsection{Construction of reduced forms}
We give some basics about reduced forms.

\smallskip
  \noindent {\bf Lemma 3.1.}
  {\em  Suppose that $w \ne \infty$.
  The reduced form of $w$  is unique.
   The reduced form of $\infty$ is either $-1/0$ or $1/0$.
}

\smallskip
\noindent {\em Proof.}
Let $a/b$ and $a'/b'$ be the reduced forms of $w$. It is clear that $ab'=a'b$.
 Following our definition of reduced form, $G_q$ contains the following two elements,
    {\tiny $\left (
\begin{array}{rr}
a & x \\
b & y \\
\end{array}
\right )$} and
{\tiny $\left (
\begin{array}{rr}
a' & u \\
b' & v \\
\end{array}
\right )$}.
An easy calculation gives
 {\tiny $\left (
\begin{array}{rr}
a & x \\
b & y \\
\end{array}
\right )^{-1}$}
{\tiny $\left (
\begin{array}{rr}
a' & u \\
b' & v \\
\end{array}
\right )$}
 $= $
  {\tiny $\left (
\begin{array}{rr}
a'y-b'x & yu-xv \\
0 & av-bu \\
\end{array}
\right )$} $\in G_q$.
This element fixes $\infty$. Since the stabiliser of $\infty$  in $G_q$ is
 generated by $T$ (see (1.1)), we have the following.
  $$\left (
\begin{array}{rr}
a & x \\
b & y \\
\end{array}
\right )^{-1}
\left (
\begin{array}{rr}
a' & u \\
b' & v \\
\end{array}
\right )
=
\left (\begin{array}{rc}
1 & m\lambda_q \\
0 & 1 \\
\end{array}
\right ).\eqno (3.1)$$

\smallskip
\noindent
 It follows from (3.1) that $a= a'$ and that $b= b'$.
 This completes the proof of the lemma.\qed

  \smallskip
  \noindent {\bf Lemma 3.2.}
    {\em Let $P$ be an ideal $q$-gon with cusps $\{c_1, c_2, \cdots , c_q\}$
     arranged in increasing order
 $($if $\infty$ is a cusp, then $c_1 = -\infty = -1/0$ if $P$ lies
   in the left half plane and $c_q=\infty=1/0$ if $P$ lies in the right half plane$)$.
   Let $a_i/b_i$ be the  reduced form of $c_i$ and $\frac{a_0}{b_0}= \frac{-a_q}{-b_q}$.
   Then}
  $$ a_i= \lambda_q a_{i-1} -a_{i-2} \mbox{ and }
                          b_i = \lambda_q  b_{i-1}-b_{i-2} \mbox{ for } 2\le i \le q-1.
                          \eqno(3.2)$$

\noindent
{\em Proof.}
Without loss of generality, we  assume that $P$
lies in the right half plane.
   Let
                                 $A = T^{-1}S= $
                               {\tiny $\left (\begin{array}{cc}
                            \lambda_q & 1     \\
                            -1\,\,\,\,\, & 0\\
                           \end{array}\right ) $} $ \in G_q$  and
                      $ w_i = {\tiny \left (\begin{array}{cc}
                            a_{i-1} & a_{i-2}     \\
                            b_{i-1} & b_{i-2}\\
                             \end{array}\right )},\, 2\le i \le q$.
                             Note that
                            $Ae^{(q-1)\pi i/q} = e^{(q-1)\pi i/q}$
and that $A$ is a counter-clockwise rotation about   $ e^{(q-1)\pi i/q}$
 of angle $2\pi/q$.

 \smallskip
(A) Suppose that the depth of $P$ is 1. Then
  $(c_1, c_q)$ is $(0, \infty)$ and $e^{\pi i /q}$ is the odd vertex of $P$.
     We apply   mathematical induction
                              as follows.

 \begin{enumerate}
 \item[(i)] It is clear that
  $w_2 = S^{-1}\in G_q$,
                               $w_2 Aw_{2}^{-1} $ fixes $    e^{\pi i/q}$  and that
                                 $w_3 = w_2 A \in G_q$.

 \smallskip
 \item[(ii)]    Suppose that  $w_{i-1}\in G_q$, $w_{i-1} A w_{i-1}^{-1} $ fixes $e^{\pi i/q}$,  and that $w_{i} = w_{i-1}A\in G_q$.

  \end{enumerate}

 \noindent    By (ii),
    $w_i Aw_i^{-1} = w_{i-1} A w_{i-1}^{-1}$ fixes  $    e^{\pi i/q}$ and   $w_i Aw_i^{-1}$
    is a counter-clockwise rotation about  $  e^{\pi i/q}$ of angle $2\pi/q$.
      Hence
    $w_{i} A w_{i}^{-1}$ sends  $(c_{i-1}, c_{i-2}) =
    (a_{i-1}/b_{i-1},  a_{i-2}/b_{i-2} ) $
     to  $(c_{i}, c_{i-1})   =  (a_{i}/b_{i},  a_{i-1}/b_{i-1} )$
  and
 $c_i$ is given by the first column of
 the following matrix      $$ w_iA w_i^{-1} w_i=
     \left (\begin{array}{cc}
                            a_{i-1} & a_{i-2}     \\
                            b_{i-1} & b_{i-2}\\
                             \end{array}\right )A
                             =
     \left (\begin{array}{cc}
                           \lambda_q a_{i-1} -a_{i-2} & a_{i-1}     \\
                           \lambda_q  b_{i-1}-b_{i-2} & b_{i-1}\\
                             \end{array}\right )\in G_q.\eqno(3.3)$$
 To be more precise,
$ (\lambda_q a_{i-1} -a_{i-2})/(
                          \lambda_q  b_{i-1}-b_{i-2}) = a_i/b_i = c_i$.
  Since
  $c_{i} >c_{i-1}$, det$\, w_i A>0$,
  and  $a_k, b_k \ge 0 $ for all $k$,
     one has    $\lambda_q  b_{i-1}-b_{i-2} >0$.
    Hence
   $ (\lambda_q a_{i-1} -a_{i-2})/(
                          \lambda_q  b_{i-1}-b_{i-2}) $ is in reduced form
                           (subsection 3.1).
                            Since the reduced form of $c_i$ is unique,   one has
                                   $ a_i= \lambda_q a_{i-1} -a_{i-2} $,
                          $b_i = \lambda_q  b_{i-1}-b_{i-2} $
                          and
                          $w_{i+1} =w_i A\in G_q$.
          By induction, we conclude that if the depth of $P$ is one,
 then  $ a_i= \lambda_q a_{i-1} -a_{i-2} $ and
                          $b_i = \lambda_q  b_{i-1}-b_{i-2} $.

\smallskip
 (B) Suppose  that the depth of $P$ is 2 or more.
                     It follows that    $(c_1, c_q)$ is  a side of a $q$-gon
  of depth $d(P)-1$.   By induction on the depth of $P$, we may assume that
   the reduced forms of $c_1 = a_1/b_1$ and $c_q= a_q/b_q$ are known and that
                     $w_2 =  {\tiny \left (\begin{array}{cc}
                            a_{1} & -a_{q}     \\
                            b_{1} & -b_{q}\\
                             \end{array}\right )}$ $\in G_q$.
                             Note that $v_0 = w_2 e^{(q-1)\pi i/q}$ is the odd vertex of $P$
                              and that
                            $w_2Aw_{2}^{-1} $ fixes $ v_0$.
                            We may now apply induction as above (see case (A)) and conclude that
                           $ a_i= \lambda_q a_{i-1} -a_{i-2} $ and that
                          $b_i = \lambda_q  b_{i-1}-b_{i-2} $.
                           This completes the proof of the lemma.\qed

\smallskip
 Lemma 3.2 tells us how to generate the reduced forms of all positive
  cusps starting with $(0/1, 1/0)$ and
 all negative
  cusps starting with $( -1/0, 0/1)$. Namely,
  $$   \frac{a_i}{b_i} = \frac{\lambda_q a_{i-1}-a_{i-2}}
  {\lambda_q b_{i-1}-b_{i-2}}\,,\,\,\, i = 2,3, \cdots, q-1.\eqno(3.4)$$

\smallskip
\noindent
Equation (3.4) allows us
 to write the $c_i$'s $(1\le i \le q$) in terms of $c_1$ and $c_q$ which
 generalises the construction of the Farey sequence. In the
 case $q=6$,  $ \lambda_6 = \sqrt 3, $  the reduced forms of the $c_i$'s $(1\le i \le 6$) in terms of $c_1= a_1/b_1 $ and $c_6= a_6/b_6$ are given by
 $$\frac{a_1}{b_1},\,\,\,\,
 \frac{  \sqrt 3 a_1+  a_6}{  \sqrt 3 b_1+  b_6},\,\,\,\,
 \frac{2 a_1+  \sqrt 3   a_6}{2 b_1+  \sqrt 3  b_6},\,\,\,\,
 \frac{  \sqrt 3  a_1+ 2 a_6}{  \sqrt 3  b_1+  2b_6},\,\,\,\,
 \frac{ a_1+  \sqrt 3   a_6}{ b_1+  \sqrt 3  b_6},\,\,\,\,
\frac{a_6}{b_6}.
 \eqno(3.5)$$

\smallskip

\noindent
 The  continuous solution of the
 equation
$b_i = \lambda_qb_{i-1} + b_{i-2} $  $( b_i \ge 0,\, 2\le i \le q-1)$
 is
$$
f(x) = \frac{b_1+b_q cos (\pi/q)}{sin (\pi/q)}
sin (\pi x/q) - b_q cos (\pi x/q)
 = A sin (\pi  x /q + w),
 \,\,f(i) = b_i,
 \eqno(3.6)$$

\noindent for some $A$ and $ w$.
Note that $f(x)$ is a periodic function of period $2q$,  $f(0) = -b_q <0$, and $f(1) = b_1 >0$.
Further,  $f(x)$ is concave down
 in the interval $[0, q]$.  As a consequence,
$$ b_i\ge b_1\mbox { and } b_i \ge b_q \mbox{ for all }  i.\eqno(3.7)$$


   \smallskip
   \noindent
 {\bf Lemma 3.3.} {\em  Let $a/b$  be in reduced form. Suppose that  $ab \ne  0$.
  Then $ |a|
 \ge 1$ and $ b\ge 1$.}

 \smallskip
 \noindent {\em Proof.}
  Since $a/b$ is in reduced form, $a/b = x_i $ is a cusp  of some $q$-gon $P$.
  Let $\{x_1, x_2, \cdots, x_q\}$ be the set of cusps of $P$.
    We  prove that $b \ge 1$ by induction on the depth of $P$.
    In the case  the depth of $P$ is 1,
 the set of cusps of $P$ is either $\{0/1, 1/\lambda_q,\cdots, 0/1\}$
  if $P$ lies in the right half plane or
   $\{-0/1, \cdots, -1/\lambda_q, 0/1\}$ if $P$ lies in the left half plane.
  By (3.7),  $b\ge 1$.   Suppose that our assertion holds when the depth of $P$ is $n$.
   We now consider the case $d(P) = n+1$ and that
    $a/b$ is a cusp of $P$.
      In the case that $a/b$ is either $x_1$ or $x_q$,
      $a/b$ is also a cusp of a $q$-gon of depth $n$.
       By inductive hypothesis, $b \ge 1$.
        Hence we shall assume
         that $a/b = x_i$, where $2\le i\le q-1$.
          By (3.7), $b$ is larger than the denominators
           of the reduced forms of $x_1$ and $x_q$.
            As $x_1$ and $x_q$ are cusps of a $q$-gon of depth $n$,
             their denominators are at least 1 by inductive hypothesis.
              Hence $b \ge 1$.
              Since $a/b$ is in reduced form, $ b/|a|$
               is also in reduced form (see subsection 3.1). Repeat the above argument, one
                has $|a|\ge 1$.
\qed

  \smallskip

\subsection{ Pseudo Euclidean algorithm}
Let $x, y  \in \Bbb Z[\lambda_q]$. Suppose that $y\ne 0$. There exists a unique
 integer  $m  \in \Bbb Z$ such that $x = y(m\lambda_q) + r $,
  $-|y\lambda_q|/2 < r\le |y\lambda_q|/2$.
 We  call such an algorithm  {\em pseudo Euclidean} (PEA).
  Let $a, b \in \Bbb Z[\lambda_q]$, where  $a b\ne  0$.
   Apply the pseudo Euclidean algorithm repeatedly,
$$ a = b(m_0\lambda_q) + r_1,$$
    $$ b = r_1(m_1\lambda_q) + r_2, $$
    $$\vdots\eqno(3.8)$$
  $$  r_{k-1} = r_k(m_k\lambda_q) + r_{k+1}. $$

 \smallskip
   \noindent
  Let $r_0=b$.  If the (PEA) terminates,
    that is $r_n \ne 0$, $r_{n+1} =0$ for some $n$, we define
    $$(a, b)_q = |r_n|.\eqno(3.9)$$

   \smallskip
    \noindent   If the (PEA) does not terminate, we define
       $(a, b)_q = 0.$
       Define further that $(a,0)_q = (0, a)_q= |a|$.
       One sees easily that
       (i) $(a,b)_q= (b,a)_q$,  (ii) $(a,b)_q = (-a, -b)_q = (-a,b)_q=(a,-b)_q$.

    \smallskip

   \noindent
 {\bf Lemma  3.4.} {\em Let $a/b$ be in reduced form. If $b=1$, then
   $a=m\lambda_q$ for some $m\in \Bbb Z$.}

   \smallskip
   \noindent  {\em Proof.}
    Suppose that $a\ne 0$. Since $a/b$ is in reduced form,
   the (PEA)
     (see (3.8)) implies that
       $r_1/b$  is  in reduced form and that
       $|r_1| < b\lambda_q/2 < b$. Since $b=1$, one has $|r_1|<1$.
        By Lemma 3.3, $r_1=0$.  Hence $a=m\lambda_q$ for some $m\in \Bbb Z$.
        This completes the proof of the lemma.\qed

  \smallskip
\noindent {\bf Remark.} The converse of Lemma 3.4 is not true
as $G_4$ and $G_6$ possess infinitely many reduced forms $m\lambda_q/b$,
    where $m, b\in \Bbb N$, $b >1$.

 \smallskip
 \noindent {\bf Lemma 3.5.} {\em Let
  $ c/d$ and $ a/b$ be the reduced forms of $x$ and $y$ respectively, where
   $x <y$. Then $g=$
  {\tiny $\left (
\begin{array}{rr}
a & c \\
b & d \\
\end{array}
\right )$} $\in G_q$ if and only if $(x,y)$ is an even line if and only if
 $ad-bc=1$.}

\smallskip
\noindent {\em Proof.} (i) Suppose that $g \in G_q$. Then $g(0, \infty)$ is an even line.
 Equivalently, $(x,y)$ is an even line. (ii) Suppose that $(x, y)$ is an even line.
  Then $ A(0,\infty) = (x, y)$ for some $A\in G_q$. An easy study of Lemma 3.1 and  the
   matrix form of $ A(0,\infty) = (x, y)$   shows that
    $ad-bc = 1$. (iii) Suppose that $ad-bc = 1$.
     Since $a/b$ and $c/d$ are in reduced forms,
     $G_q$ contains elements of the following forms,
     {\tiny $\left (
\begin{array}{rr}
a & x \\
b & y \\
\end{array}
\right )$} and
{\tiny $\left (
\begin{array}{rr}
c & u \\
d & v \\
\end{array}
\right )$}. An easy calculation show that
 {\tiny $\left (
\begin{array}{rr}
a & x \\
b & y \\
\end{array}
\right )^{-1}$}
{\tiny $\left (
\begin{array}{rr}
c & u \\
d & v \\
\end{array}
\right )$}
  $ = $
{\tiny $\left (
\begin{array}{cr}
cy-dx & * \\
1 & *\\
\end{array}
\right )$}. By Lemma 3.4, one has $cy-dx = m\lambda_q$ for some $m \in \Bbb Z$.
 It follows  that
 $ {\tiny \left (
\begin{array}{rr}
a & x \\
b & y \\
\end{array}
\right )^{-1}
g =
  \left (
\begin{array}{rr}
a & x \\
b & y \\
\end{array}
\right )^{-1}
\left (
\begin{array}{rr}
a & c\\
b & d \\
\end{array}
\right ) =
\left (
\begin{array}{rc}
1 & m\lambda_q \\
0 & 1\\
\end{array}
\right )}\in G_q.$ Hence $g\in G_q$.\qed

\smallskip
\noindent {\bf Lemma 3.6.} {\em Suppose that $(a,b)_q =r\ne 0$. Then
there exists some $g\in G_q$ such that
$g
{\tiny \left ( \begin{array}{c}
a\\
b\\\end{array}
\right)=
\left ( \begin{array}{c}
r\\
0\\
\end{array}\right )}$.
In particular, $(a/r, b/r)_q =1$ and $\frac{a/r} {|b/r|}$ is in reduced form.
}
\smallskip

\noindent {\em Proof.}
 The lemma follows  from (3.8), (3.9)  and the observation that the matrix form of the equation $x=y(m\lambda_q)+z$ is
{\tiny $
\left ( \begin{array}{cc}
1 & -m\lambda_q\\
 0 & 1\\
\end{array}
\right )
\left (\begin{array}{c}
x\\
y\\
\end{array}
\right ) =
\left ( \begin{array}{c}
z\\
y\\
\end{array}
\right )$}, where
{\tiny $
\left ( \begin{array}{cc}
1 & -m\lambda_q\\
 0 & 1\\
\end{array}
\right )
$} $\in G_q$.
\qed

\smallskip

    \noindent {\bf Proposition 3.7.} {\em Let
   $A =$  {\tiny $\left (
\begin{array}{cc}
a&c\\
b &d \\
\end{array}
\right )$ } $\in PSL(2, \Bbb Z[\lambda_q])$.
 Then $A\in G_q$ if and only if
  $(a, b)_q= (c,d)_q=1$. In particular,
   if $(a, b)_q <1$  or $|b|<1$, then $A$ is not in $G_q$.}
\smallskip

\noindent {\em Proof.}
Suppose that $A \in G_q$. Then
 $a/|b|$ is in reduced form (see subsection 3.1). By (3.8),
  $r_{n}/|r_{n+1}|$ is in reduced form for every  $n$. By Lemma 3.3, $|r_n|\ge 1$ and
   $|r_{n+1}|\ge 1$ whenever $r_nr_{n+1} \ne 0$.
Since $a$ and $b$ are finite and $\lambda_q/2 <1$, an easy observation of (3.8)
 implies that there exists some $m $ such that $|r_k| <1$ whenever $k \ge m$.
 Hence there exists a $d$ such that $r_{d+1} =0$ and that $r_d\ne 0$.
  By Lemma 3.1, $r_d =\pm 1$.
   Equivalently, $(a, b)_q = 1$.
    Since $A\in G_q$, $A^{-1}S \in G_q$. Since the transpose of $S$ and $T$ are
     members of $G_q$, the transpose of $A^{-1}S$ is also an element of $G_q$,
      the first column of  the transpose of $A^{-1}S$ is
{\tiny $\left (
\begin{array}{c}
c\\
d \\
\end{array}
\right )$}. This implies that $c/|d|$ is in reduced form. Similar to the above,
 one can show that  $(c, d)_q =1$.

\smallskip
Conversely, suppose that $(a,b)_q= (c,d)_q = 1$.
 Replace $A$ by $-A$ if necessary, we may assume that $d\ge 0$.
  Replace $A$ by $wAw^{-1}$ if necessary
 (see subsection 3.1 for $w$), we may assume that $b\ge 0$ and that $d\ge0$. By Lemma 3.6,   both $a/b$ and $c/d$ are in reduced forms.
  By Lemma 3.5, we have $A \in G_q$.
  This completes the proof of the lemma.\qed

 \smallskip
 \noindent
  {\bf Example 3.8.} Let $q=5 $ and let $\lambda = \lambda_5$.
   Note that $\lambda^2= \lambda +1$.
 {\tiny$\left (
\begin{array}{cc}
4\lambda -1&\lambda +1\\
 3&\lambda \\
\end{array}
\right )$} is not an element of $G_5$ as
 $(4\lambda-1, 3)_5 = \lambda-1<1$.

   \section{Hecke-Farey Symbols and special polygons}
\subsection{ $r$-clusters}
 Let $\Phi$ be the hyperbolic triangle
 with vertices $0, e^{\pi i/q}$ and $\infty$. $\Phi$ is a fundamental domain of $G_q$.
 The $G_q$ translates of $\Phi$ are called {\em special triangles}.
  For each  divisor $r$ of $q$ $(1\le r< q)$,
 set
 $\Phi_r = \Phi \cup R\Phi \cup R^2\Phi \cup \cdots \cup R^{r-1}\Phi.$
$\Phi_r$ is a union of $r$ copies of special triangles. These  $r$ special triangles meet
 at the odd vertex $e^{\pi i/q}$.
 The $G_q$ translates of $\Phi_r$
 are called the {\em $r$-clusters} (a 1-cluster is a special triangle).
 Let $\Delta_r$
  be an $r$-cluster. It is clear that
  (i) $\Delta_r$ has $(r+1)$ cusps and one odd vertex $y$, (ii)
  the boundary of $\Delta_r$ has $r$ even lines and  two odd edges,
 (iii)  the two odd edges   of $\Delta_r$ meet each other at $y$
  with
   vertex angle $2r\pi/q$. The cusps are called the {\em free vertices} of $\Delta_r$.

   Let $\{ v_1 = 0/1 , v_2, \cdots, v_{q-1} ,  v_q= 1/0\}$
  be the set of cusps  of the depth 1 $q$-gon
  in the right half plane
   where the $v_i$'s are arranged in increasing order.
   Following the definition of $\Phi_r$,
    the set of free vertices (cusps) of  $\Phi_r$ is
    $$\Phi_r \,:\,\{ 0= v_1 , v_2, \cdots, v_{r}, v_q=\infty \}.\eqno(4.1a)$$
     The  odd vertex of $\Phi_r $ is $ e^{\pi i/q}$.
      The odd edges are $(e^{\pi i/q}, v_{r})
     $ and $(e^{\pi i/q}, v_q)=(e^{\pi i/q}, \infty)$.
      Note that $\Phi_r$  is convex as it is  the convex hull of
       $\{ v_1 , v_2, \cdots, v_{r},  e^{\pi i /q}, v_q \}$.
 Let  $\Psi_r$ be the mirror image of $\Phi_r$ (with respect to the
  $y$-axis).
 The set of free vertices of $\Psi$ is
  $$ \Psi_r \,:\,\{-\infty =  -v_q , -v_{r}, \cdots , -v_2, -v_1 = 0\}.\eqno(4.1b)$$

 \noindent The odd vertex of $\Psi_r$ is
   $y= e^{(q-1)\pi i/q}$ and the  odd edges of $\Psi_r$  are $(y, -v_{r})
     $ and  $(y, -v_q)$.  Let $A$ be given as in Lemma 3.2.
   It is clear that   $\Psi_r = A^{1-r}S\Phi$.

   \smallskip

      \noindent {\bf Example 4.1.}
Let $q=6$, $\lambda= \lambda_6 =\sqrt 3$.
See {\bf  Figure  1a} for $\Phi_2$ and {\bf Figure 1b} for $\Psi_2$.

  \allowdisplaybreaks

\begin{center}
\begin{figure}
\beginpicture

\setcoordinatesystem units <9pt,9pt>
\setplotarea x from  -8 to 63, y from -1 to 9

\setlinear \plot 0  0    0 8 /

\setlinear \plot 8.66 5  8.66 8 /

\circulararc -120 degrees from 0 0 center at 5.77 0

\circulararc -180 degrees from 0 0 center at 2.885 0

\circulararc -60 degrees from 5.77 0 center at 11.55 0


\put { \small    $0$} at   0 -2
\put { \small $1/\lambda$} at   5.77 -2

\put {\small$ \,e^{\pi /6}$} at   10  5


\setlinear \plot 20  5    20 8 /

\setlinear \plot 28.66 0   28.66 8 /
\put {\small $ e^{5\pi /6}$} at   18  5

\circulararc 120 degrees from 28.66 0 center at 22.89 0

\circulararc 180 degrees from 28.66 0 center at 25.775 0
\put {\small $ -1/\lambda$} at  22.49 -2

\circulararc 60 degrees from 22.89 0 center at 17.11 0

\put { \small $0$} at   28.66 -2

\put {\bf Figure 1a} at 4.77  -5
\put {\bf Figure 1b} at  24.77 -5
\endpicture
\end{figure}
\end{center}

\vspace{-.5cm}

\subsection{Special polygons}
A convex hyperbolic polygon $P$  of $\Bbb H$ is
 a union of some $q$-gons  and a finite number of  $r_i$-clusters
 ($r_i|q$, $1\le r_i <q$).
  The $q$-gons and  the $r_i$-clusters of $P$ are
     called the
    {\em tiles}.
     The tiles intersect each other (if any) at either free vertices (cusps) or
     even lines.
A {\em special polygon} $M_X= (P, I_X)$ of $\Bbb H$ is a convex hyperbolic polygon $P$
   together with a set of
    side pairings $I_X$ satisfying the rules below.

    \begin{enumerate}
    \item[(S1)] An odd edge $e$  is always paired with an odd edge $f$
     (in the same $r$-cluster)
    and
    makes an internal angle  $2r\pi/q$ with $f$.
      The vertex  where
       $e$ and $f$  meet is  an  {\em odd vertex} of $P$.
      Both $e$ and $f$ are considered as sides of $P$, and are called its {\em odd
       sides}.

     \item[(S2)] Let $e$ and $f$ be two even edges in the boundary of $P$ forming an
      even line. Then either (i) $e$ is paired with $f$,
       both $e$ and $f$ are considered as sides of $P$, and  are
       called its {\em even sides},
        the point where $e$ and $f$  meet is
        an  {\em even  vertex} of $P$,
        or  (ii) $e$ and $f$ form a {\em free side}
       of $P$, and this free side is paired with another free side of $P$.
\item[(S3)] $0$ and $\infty$ are vertices of $P$.
\end{enumerate}

\noindent Let $M_X=(P, I_X) $ be a special polygon. The cusps in $P$ are called the {\em free vertices}
 of $M_X$.

\subsection{ Hecke-Farey sequences and symbols.}
    A  {\em Hecke-Farey sequence}   is a  finite sequence of cyclically arranged
  numbers in increasing order
 $\{-\infty, x_0, x_1, \cdots, x_n, \infty\}$
  such that
  \begin{enumerate}
  \item[(a)] $ x_i \in \Bbb Q [\lambda_q]=\Bbb Q(\lambda_q)$, $x_i=0$ for some $i$, $0\le i \le n$,
  \item[(b)]   $x_i = a_i/b_i$ is in reduced form  for every $i$ ($x_{-1} =-\infty
   = -1/0$ and $x _{n+1} = \infty = 1/0$),
 \item[(c)] if $a_{i+1}b_i -a_ib_{i+1} \ne 1$, then
  there exists an element $g\in G_q$ and an $r$-cluster $\Phi_r$  $(r|q,\,1 < r<q)$ such that
   $g(v_{r}) = x_i, \,g(v_q) = x_{i+1}$ (see (4.1a) for $v_{r} $ and $v_q$).
\item[(d)] $\{x_i, x_{i+1}\}$ is called an {\em interval}.
 Intervals described as in (c) are called
  {\em r-intervals}.
  $\{x_i, x_{i+1}\}$ is called an {\em ordinary interval}
if $a_{i+1}b_i -a_ib_{i+1} = 1$.
  By Lemma 3.5, if $\{x_i, x_{i+1}\}$ is an ordinary interval,
    then
 the
   hyperbolic
    line $(x_i, x_{i+1})$ is an even line.

  \end{enumerate}
One can show easily that  that the  $g$'s and $r$'s in (c) of the above are
 unique
 ({the only element of $G_q$ that fixes two or more  points of $\Bbb H$ is the identity}).
A  {\em Hecke-Farey symbol} is a Hecke-Farey sequence together with
 an additional {\em side pairing}
 on each consecutive pair of $x_i$'s. {\em To avoid triviality, we insist that
  a Hecke-Farey symbol
   must have at least two distinct  side pairings} (see Discussion 4.5).
  In the case $a_{i+1}b_i -a_ib_{i+1} \ne 1$,
 the additional side pairing   of the $r$-interval  $\{x_i, x_{i+1}\}$ is denoted by
 $$ {x_i}\,\, _{_{_{\smile}} }   \hspace{-.4cm} _{_{_{_{_{_{_{e_{r(g)}}}}}}}}  x_{i+1},
 \eqno(4.2) $$
where $g$  and $r$ are  given as in (c) of the above.
 In the case $a_{i+1}b_i -a_ib_{i+1} =1$,
  the  additional side pairing
 of the interval   $\{x_i,x_{i+1}\}$
 can be  either one of the following three types :
  $$ {x_i} _{_{_{\smile}}} \  \hspace{-.37cm} _{_{_{_{_{_{_{\circ}}}}}}}
  x_{i+1}, \,\,
{x_i} _{_{_{\smile}}} \  \hspace{-.37cm} _{_{_{_{_{_{_{\bullet}}}}}}}
 x_{i+1} ,\,\,
  {x_i} _{_{_{\smile}} }\  \hspace{-.37cm} _{_{_{_{_{_{_{a}}}}}}}  x_{i+1},\,\,
\eqno(4.3) $$
 where $a$ is a natural  number.
 Each natural number $a$ occurs
 exactly twice or not at all.
 The actual values of the $a$'s are unimportant, it is
 the pairing induced on the consecutive pairs that matters.
We shall now give a detailed description of these side pairings.

\begin{enumerate}

\item[(i)] Let $\{x_i, x_{i+1}\}$  be an $r$-interval and let $g, \Phi_r$ be given as in
 (c) of the above. By (c),
 one has $g(v_{r}) = x_i$ and  $g(v_q)= x_{i+1}$.
 Set $g(e^{\pi i/q}) = y$.
  $$ {x_i}\,\, _{_{_{\smile}} }   \hspace{-.4cm} _{_{_{_{_{_{_{e_{r(g)}}}}}}}}  x_{i+1}
  \eqno(4.4)$$
  is the side pairing that
 fixes $y$ and sends the odd edge $(y, x_{i+1})$ to $(y,x_i)$.
  Since $R$ is a counterclockwise rotation about $e^{\pi i/q}$
  of angle $2\pi/q$ and $(y ,x_i)$ and $ (y, x_{i+1})$ form an internal angle $2r\pi/q$,
  this side pairing must be  $gR^rg^{-1} \in G_q$
  ({the only element of $G_q$ that fixes two or more  points of $\Bbb H$ is the identity}).

\smallskip

\item[(ii)]
Let $\{x_i, x_{i+1}\}$ be an ordinary interval.
By (d) of the above,  $(x_i, x_{i+1})$ is an even line.
 Hence there exists a unique
  $A\in G_q$ such that $A(0)= x_i$ and $A(\infty) = x_{i+1}$.
   Set
    $A(e^{\pi i/q})=y$.  $(y, x_{i+1})$ and $(y, x_i)$ are  odd edges.
  $$ {x_i} _{_{_{\smile}}} \  \hspace{-.37cm} _{_{_{_{_{_{_{\bullet}}}}}}}
  x_{i+1}\eqno (4.5)$$
is the side pairing that fixes $y$ and sends the odd edge $(y,x_{i+1})$ to $(y,x_i)$.
 It is clear that this side pairing must be
    $$ h= ARA^{-1}=
\left (
\begin{array}{rr}
a_{i+1} & a_i \\
b_{i+1} & b_i \\
\end{array}
\right )
  \left (
\begin{array}{rc}
0 & 1\\
-1 & \lambda_q \\
\end{array}
\right )
\left (
\begin{array}{rr}
a_{i+1} & a_i \\
b_{i+1} & b_i \\
\end{array}
\right ) ^{-1}
\in G_q \,.\eqno(4.5a)
$$

\smallskip

\item[(iii)]
Let $\{x_i, x_{i+1}\}$ be an ordinary interval.
 By (d) of the above,  $(x_i, x_{i+1})$ is an even line.
 Hence there exists a unique
  $A\in G_q$ such that $A(0)= x_i$ and $A(\infty) = x_{i+1}$.
   Set
    $A(\sqrt {-1} )=y$.  $(y, x_{i+1})$ and $(y, x_i)$ are even edges.
  $$ {x_i} _{_{_{\smile}}} \  \hspace{-.37cm} _{_{_{_{_{_{_{\circ}}}}}}}
  x_{i+1}\eqno (4.6)$$
is the side pairing that fixes $y$ and sends the even edge $(y,x_{i+1})$ to $(y,x_i)$.
 It is clear that this side pairing must be
 $$ w= ASA^{-1}=
\left (
\begin{array}{rr}
a_{i+1} & a_i \\
b_{i+1} & b_i \\
\end{array}
\right )
  \left (
\begin{array}{rr}
0 & 1\\
-1 & 0 \\
\end{array}
\right )
\left (
\begin{array}{rr}
a_{i+1} & a_i \\
b_{i+1} & b_i \\
\end{array}
\right )
^{-1}
\in G_q \,.\eqno(4.6a)
$$

\smallskip
\item[(iv)]
Let $\{x_i, x_{i+1}\}$ and
 $\{x_j, x_{j+1}\}$ be two ordinary intervals.
  By (d) of the above,
  $(x_j, x_{j+1})$ and  $(x_i, x_{i+1})$ are even lines.
  $$ {x_i} _{_{_{\smile}} }\  \hspace{-.37cm} _{_{_{_{_{_{_{a}}}}}}}  x_{i+1}\,,\,\,\,
  {x_j} _{_{_{\smile}} }\  \hspace{-.37cm} _{_{_{_{_{_{_{_a}}}}}}}  x_{j+1}\eqno(4.7)$$
is the side pairing that sends $x_{j+1}$ to $x_i$ and $x_j$ to $x_{i+1}$
 (equivalently, pairs the even lines $(x_j, x_{j+1})$ and  $(x_i, x_{i+1})$).
  It is clear that this side pairing must be
$$k= \left (
\begin{array}{rr}
a_j & -a_{j+1} \\
b_j & -b_{j+1} \\
\end{array}
\right )
  \left (
\begin{array}{rr}
a_{i+1} & a_i\\
b_{i+1} & b_i \\
\end{array}
\right )^{-1}\in G_q .\eqno(4.7a)
$$

\end{enumerate}

\noindent \noindent
 {\bf Definition 4.2.}
  The intervals in (ii)  are called {\em odd intervals}, the intervals in
  (iii) are called {\em even intervals},
and the  intervals in (iv) are called {\em free intervals.}

     \smallskip
     \noindent
     {\bf Example 4.3.} Let $q=6, \,\lambda = \lambda_6=  \sqrt 3$,
      $g = $
 ${\tiny \left (
\begin{array}{rr}
1 & 0 \\
\lambda  & 1 \\
\end{array}
\right )}$. The following is a Hecke-Farey symbol
 (see Figure 2).

$$ M_X=
  \{
   {-\infty }\,
   _{_{_{\smile}} }\  \hspace{-.37cm}_{ _{_{_{_{_{_{_{\circ}}}}}}} }
 \,\,
   {0/1 }\,
   _{_{_{\smile}} }\  \hspace{-.37cm}_{ _{_{_{_{_{_{_{\circ}}}}}}} }
 {1/2\lambda}\,\,
_{_{_{\smile}} }\  \hspace{-.5cm}_{ _{_{_{_{_{_{_{e_2(g)}}}}}}}}
1/\lambda\,
 _{_{_{\smile}}}\  \hspace{-.36cm} _{_{_{_{_{_{_{_{1}}}}}}}}
\,\lambda/2\,
 _{_{_{\smile}}}\  \hspace{-.36cm} _{_{_{_{_{_{_{_{1}}}}}}}}
\,2/\lambda \,
_{_{_{\smile}} }\  \hspace{-.37cm}_{ _{_{_{_{_{_{_{\bullet}}}}}}}}
\,\lambda/1\,
_{_{_{\smile}} }\  \hspace{-.37cm}_{ _{_{_{_{_{_{_{\circ}}}}}}}}
{\infty}
 \} .\eqno(4.8)$$

\noindent {\em Proof.} The set of cusps of $\Phi_2$ is $\{0/1, 1/\lambda,
\infty\}$.
One must verify  that the interval $ \{1/2\lambda, 1/\lambda\}$
 satisfies (c) of the above.
   This  can be checked  easily as
    $g(1/\lambda) = 1/2\lambda$ and $g(\infty) = 1/\lambda$.
    Note that $g\Phi_2$ is the 2-cluster with odd vertex $y$.
      \qed

\allowdisplaybreaks

\begin{center}
\begin{figure}
\beginpicture

\setcoordinatesystem units <4pt,4pt>

\setplotarea x from  -13 to 83, y from -1 to 20

\setlinear \plot 0 0 0 18 /
\setlinear \plot 75 0 75 18 /

\setlinear \plot 0 0 75  0  /

\circulararc -180 degrees from 0 0 center at 12.5   0
\circulararc -180 degrees from 25 0 center at 31.25    0
\circulararc 180 degrees from 50 0 center at  43.75    0
\circulararc 180 degrees from 75 0 center at 62.5   0



\circulararc -157 degrees from 00 0 center at 8.33 0
\circulararc -180 degrees from 00 0 center at 6.25 0

\circulararc -96 degrees from 12.5 0 center at 15.626 0

\circulararc 142 degrees from 25 0 center at 20 0

\circulararc 158 degrees from 75 0 center at 66.67  0
\circulararc  -141 degrees from 50 0 center at 55  0

\put {\tiny $ 0/1$} at 0 -3

\put {\tiny $ 1/\lambda$} at 25 -3
\put {\tiny $ 1/2\lambda$} at 12.5 -3
\put {\tiny $ y$}         at 17    1.8
\put {\tiny $ w$}         at 60    1.8

\put {\tiny $ \lambda/2$} at 37.25 -3
\put {\tiny $ 2/\lambda$} at 50 -3
\put {\tiny $ \lambda/1$} at 75 -3

\put {\bf  Figure 2} at 37.25 -7

\endpicture
\end{figure}
\end{center}

\vspace{-.5cm}

\subsection{Special polygons and Hecke-Farey symbols}
We now relate Hecke-Farey symbols to special polygons. It is a generalisation of
 a theorem of [K3] and is proved in an analogous way
  (see Example 4.3 and Figure 2 for an example).

\smallskip

\noindent {\bf Proposition 4.4.} {\em There is a one to one correspondence between
 the set of special  polygons   and the set of Hecke-Farey symbols.}
\smallskip

\noindent {\em Proof.}
Let $M_X= (P, I_X) $ be a special polygon. The free vertices (cusps) of $P$ form
 a sequence $\{-\infty = x_{-1}, x_0, x_1,\cdots, x_n, \infty= x_{n+1}\}$.
  Let $a_k/b_k$ be the reduced form of $x_k$ and let $x_i, x_{i+1}$ be two consecutive
   terms.
   Suppose that $a_{i+1}b_i-a_ib_{i+1}\ne 1$. By Lemma 3.5, $(x_i, x_{i+1})$ is not
    an even line of $P$.
    Following the definition of convex hyperbolic polygons (see (S1)-(S3) of subsection
    4.2), one has (i)
     $x_i$ and $x_{i+1}$ are the end points of
     two odd edges $e = (y, x_i) $ and $f=(y, x_{i+1})$ of $P$ that meet at an odd vertex $y$, (ii) $e$ and $f$  are the odd edges of an $r$-cluster $g\Phi_r$
       for some $g$ and $r$. As a consequence,
        $x_i$ and $ x_{i+1}$ satisfy (c) of subsection 4.3 and gives an $r$-interval
        $\{x_i,x_{i+1}\}$.
       Hence
       the sequence
      $$ F = \{-\infty = x_{-1}, x_0, x_1,\cdots, x_n, \infty= x_{n+1}\}\eqno(4.9)$$
       is a Hecke-Farey sequence.
        In the case $\{x_i, x_{i+1}\}$ is an $r$-interval, following our  notation above and (S1) of subsection 4.2,  $e=(y, x_i) $ and $f=(y, x_{i+1})$ are paired by an element of $I_X$.
            As this side pairing must be unique, it has to be
             $gR^rg^{-1}$ (see (4.4)).
             Hence the $r$-interval must carry the side pairing  given by (4.4).
              In the case $\{x_i, x_{i+1}\}$ is an ordinary interval of  $F$,
               similar study of (S1)-(S3) of subsection 4.2 implies
                that the side pairing in $I_X$
                that does the pairing for $\{x_i, x_{i+1}\}$
                 is just the side pairing  given by (4.3). This makes $F$ into a
                  Hecke-Farey symbol.

\smallskip
Conversely,  given a Hecke-Farey  symbol, one can  construct a
    special   polygon with free vertices
    corresponding to the points of the
     Hecke-Farey sequence, odd vertices
      where the symbols $e_r(g)$ and $\bullet$
      occur and even vertices where the symbol $\circ$ occurs
       and side pairings as determined by (i)-(iv) of subsection 4.3.\qed

\smallskip
\noindent {\bf Discussion 4.5.}
Note that
$ F_1 =
  \{
   {-\infty }\,
   _{_{_{\smile}} }\  \hspace{-.37cm}_{ _{_{_{_{_{_{_{\circ}}}}}}} }
 \,\,
   {0/1 }\,
   _{_{_{\smile}} }\  \hspace{-.37cm}_{ _{_{_{_{_{_{_{\circ}}}}}}} }
{\infty}
 \} \mbox{ and }
 F_2 =
  \{
   {-\infty }\,
   _{_{_{\smile}} }\  \hspace{-.37cm}_{ _{_{_{_{_{_{_{1}}}}}}} }
 \,\,
   {0/1 }\,
   _{_{_{\smile}} }\  \hspace{-.37cm}_{ _{_{_{_{_{_{_{1}}}}}}} }
{\infty}
 \}
 $
  are not Hecke-Farey symbols since a Hecke-Farey symbol must have
   at least two side pairings (see subsection 4.3) whereas
  $F_1$ has only one side pairing $S$  and $F_2$ has only one side pairing $I_2$.
   It is clear that $F_1$ and $F_2$ do not correspond to any special polygon.
  As a consequence, Proposition 4.4 is no longer true
   if one does not insist that a Hecke-Farey symbol must carry at least two
    side pairings.

  \section{Poincar\'{e}'s polygonal theorem and independent generators}

 \noindent {\bf Theorem 5.1.}
   {\em If
     $M_X= (P, I_X)$ is a special polygon, then the set of side pairings $I_X$ generates independently
      a group $G \subseteq G_q$ such that $P$ is a fundamental domain of $G$.}

  \smallskip
  \noindent {\em Proof.}  The stabilisers of the cusps of $P$ in $G$ are generated by
   conjugates of $T^n$ for
    some $n\in \Bbb N$. Let $v\in P$ be an elliptic point of $P$. The stabiliser $ G_v
    =\left < \tau\right >$
     of $v$ in $G$ is generated by either a conjugate of $S$ or a conjugate of $R^r$
      for some $r$, where $r|q$.
        It is clear from our construction of $P$ that (i) $v$ is the intersection of
      two edges $e_1$ and $e_2$ of $P$ (see (S1) and (S2) of subsection 4.2), (ii)
       the sides $e_1$ and $e_2$ make an internal angle $\pi$ or $ 2r\pi/q$ at $v$,
        and  (iii) we may assume that $\tau$ fixes $v$ and sends
      $e_1$ to $e_2$.
       In summary,
        \begin{enumerate}
        \item[(i)] the stabiliser of a cusp of $P$ in $G$ is generated by a parabolic element,

        \item[(ii)]  $e_1$ and $e_2$ make an internal angle $\pi$ or $ 2r\pi/q$ at $v$,
         where $r|q$ and $r <q$,

        \item[(iii)] $\tau$  is the side pairing that fixes $v$ and sends
        $e_1$ to $e_2$.
        \end{enumerate}
       Hence $P$ is a Poincar\'{e} polygon. By  Poincar\'{e}'s polygonal theorem,
   $P$ is a fundamental domain of $G$
            and $I_X$ is  a set of independent generators of $G$
             (see  pp 223 of [Ma]).\qed

\smallskip

   Let $X$ be a subgroup of finite index of $G_q$.
       An  {\em admissible fundamental domain} is a special polygon
        $M_X= (P, I_X)$ such that $P$ is a fundamental domain of $X$ and $I_X$ is
         a set of independent generators of $X$.

  \smallskip
  \noindent {\bf Theorem 5.2.}
   {\em Let $X$ be a subgroup of finite index of $G_q$. Then $X$ has an admissible fundamental
    domain $M_X$.}

  \smallskip
  \noindent {\em Proof.}
The tessellation $\mathcal I^*$ (see subsection 2.1) induces a tessellation of the surface
 $X\setminus \Bbb H$.
 Let $v$ be an  odd vertex of a tile $T_0$ of  $X\setminus \Bbb H$.
 The stabiliser $X_v\subseteq X$ is  cyclic of order $d$, where $d$ is a divisor of $q$.
   In the case $d>1$, $v$ is the odd vertex of
   a  $q/d$-cluster $T_0$. Note that
  the two odd sides
     of $T_0$ must be paired by some members of $X_v$.
       In the case $d=1$, $v$ is the odd vertex of a $q$-gon
       $T_0$. In summary, the tiles of $X\setminus \Bbb H$
        are $q$-gons and $q/d$-clusters, where $d|q$.
 If $X\setminus \Bbb H$ has an  elliptic point of order 2 on the
        boundary of  $T_0$,
         then the two even edges incident to the elliptic point
         are paired by an element of $X$ and forming an edge of $X\setminus \Bbb  H$.
          {\em  It is clear that these tiles  intersect each other
          $($if any$)$ at either cusps or even lines as the two odd sides of a  $q/d$-cluster
           must be paired by some elements in $X_v$}.

\smallskip
  Finding a fundamental domain for $X$ which is a special polygon amounts to
  cutting the surface $X\setminus \Bbb H$ into its tiles ($q$-gons and $q/d$-clusters)
   and develop these tiles on $\Bbb H$ so that $0$ and $
   \infty$ are vertices of the polygon (see (S3) of subsection 4.2).

\smallskip
 \noindent (A) Let $p\,:\, \Bbb H \to X\setminus \Bbb H$ be the projection map and
 let $(0, \infty)$ be  the even line joining 0 and $\infty$.
  $p((0,\infty))$ lies on the boundary of some tile $T$ of $X\setminus \Bbb H$.
   We develop $T$ to $\Bbb H$  so that $p((0,\infty))$ is developed to $(0, \infty)$.
The other tiles of $X\setminus \Bbb H$ are then developed onto $\Bbb H$ in an inductive
 manner, piece by piece, where each new piece is adjacent to a tile that is already been
  developed.

\smallskip

\noindent (B) The determination of the first tile $P_0$.
\begin{enumerate}
\item[(i)]
 Let
$D_1 =\{\Phi_d \,:  \,\mbox  {the two odd edges  of } \Phi_d
 \mbox{ are paired by some elements of } X\}$ and
 $D_2=\{\Psi_d \,:  \,\mbox  {the two odd edges  of } \Psi_d
 \mbox{ are paired by some elements of } X\}$
  (see subsection 4.1 for $\Phi_d$ and $\Psi_d$).
  If $D_1\ne \emptyset$, let $P_0$ be the smallest $\Phi_d$  (in area) of $D_1$.
  If $D_1=\emptyset$ and $D_2\ne \emptyset$,
  let $P_0$ be the smallest $\Psi_r$  (in area) of $D_2$.
  If $D_1\cup D_2=\emptyset$,
    let $P_0$ be the depth one $q$-gon that lies in the right half plane.
    Note that 0 and $\infty$ are vertices of $P_0$ (see (S3) of subsection 4.2).
\item[(ii)]
 The reduced forms of the vertices of
  $P_0$ can be determined by Lemma 3.2. Note that the two odd edges of $P_0$
   (if any)  are
 paired by some  elements of $X$.
\end{enumerate}

\smallskip
\noindent (C)   Let $e$ be a side (even line) of $P_0$.

\begin{enumerate}
\item[(i)] If there exists another side $f$ of $P_0$ such that
 the element $g\in G_q$ which pairs $e$ and $f$ is in $X$, then
  we call  $e$ and $f$  {\em paired sides}  and add $g$ to the generating set $I_X$ of $X$.
The side pairing $g$ can be determined by  (4.7).

\item[(ii)] If the element $g$ in $G_q$ that pairs the two even edges of $e$ is in $X$,
then we call $e$ a {\em paired side} and
 $g$ is put into the generating set $I_X$ of $X$.
Note that $g$ has order 2. Such $g$ can be determined by (4.6).
\item[(iii)] Let
 $D $ be the collection of all $d$-clusters $(d<q$) attached to $e$
 and let $D_e =\{ T\in D\,:\, \mbox{the two odd edges of $T$ are paired by
  some elements of $X$}\}$.
   If  $D_e$ is not empty,
       then we  attach
     $T_0$ to $P_0$ (along $e$) to form a new polygon $P_1 = P_0\cup T_0$,
     where $T_0$ is the smallest $d$-cluster (in area) of $D_e$.
     The two odd edges of $T_0$  paired by $g \in X$ are called {\em paired sides}
      and $g$ is put into the generating set $I_X$ of $X$  as before.
       Note that $g$ can be determined by (4.4)  and (4.5).

\item[(iv)] If $e$ does not satisfy any of the above, $e$ is called an
 {\em unpaired side} of $P_0$.

\item[(v)] It is clear that only one of the four cases can occur to $e$.

\end{enumerate}

\smallskip
\noindent
(D) The inductive step.
\smallskip
\begin{enumerate}
\item[]
We have a polygon $P_n$ with some paired sides and some unpaired sides.
If there is no  unpaired sides, we are done and $P_n$ with the  side pairings is a
 special polygon of $X$. Otherwise choose any one of the unpaired sides $e$
  and adjoin a $q$-gon $T_n$ along $e$ to form a new polygon $P_{n+1} = P_n\cup T_n$.
  The reduced forms of the vertices of $T_n $ can be obtained by
   applying Lemma 3.2. $P_{n+1} $ has $q-1$ new sides (coming from $T_n$).
   We check for cases (i), (ii) and (iii) of (C) for each of the new sides.
   If after checking (i)-(iii), $P_{n+1}$ has no more unpaired sides,
    we are done. Otherwise we continue inductively.

\end{enumerate}

\smallskip

  \noindent
   Since the index of $X$  is finite, this process must stop after
    a finite number of steps. The polygon  $P$  together with the set of side pairings
     $I_X$ is a special polygon.
        By Theorem 5.1,  $(P, I_X)$ is an admissible fundamental domain of  $X$.\qed

   \smallskip


  \subsection{ Two  subgroups of  the Hecke group $G_6$ of index 3} (i)
  Let $\lambda = \lambda_6 = \sqrt 3$ and let $M_X$ be a  special polygon of $X$  given as follows (see Figure 3a).
$$
  M_X= \{
   {-\infty }_{_{_{\smile}} }\  \hspace{-.37cm}_{ _{_{_{_{_{_{_{\circ}}}}}}} }
 \,\,
   {0/1 }_{_{_{\smile}} }\  \hspace{-.37cm}_{ _{_{_{_{_{_{_{\circ}}}}}}} }
 {1/\lambda}\,\,
_{_{_{\smile}} }\  \hspace{-.32cm}_{ _{_{_{_{_{_{_{\circ}}}}}}}}
\lambda/2 \,\,
_{_{_{\smile}} }\  \hspace{-.52cm}_{ _{_{_{_{_{_{_{e_3(1)}}}}}}}}
\infty
 \}. \eqno(5.1)$$

\noindent To determine the side pairing $e_3(1)$, we follow the notation of (i)
 subsection 4.3. One has $r=3,$ $g=1$,  $x_i = \lambda/2$ and that $x_{i+1} = 1/0$.
  By (4.4), one has $e_3 (1)= R^3$. $X$ is a normal subgroup of index 3 (see Discussion 8.2
   for normality).
   The remaining side pairings can be determined easily by (4.6).
In summary,
  a set of independent generators is given by
      $$
    X = \left < S,
\left (
\begin{array}{rr}
1 & 0 \\
 \lambda & 1 \\
\end{array}
\right )
 S
\left (
\begin{array}{rr}
1 & 0\\
\lambda & 1 \\
\end{array}
\right ) ^{-1},
\left (
\begin{array}{rr}
\lambda & 1   \\
2 & \lambda \\
\end{array}
\right ) S
\left (
\begin{array}{rr}
\lambda  & 1\\
2 & \lambda \\
\end{array}
\right ) ^{-1},
\left (
\begin{array}{rr}
0 & 1\\
-1 & \lambda \\
\end{array}
\right ) ^{3}\right >.\eqno(5.2)
$$

  \allowdisplaybreaks

\begin{center}
\begin{figure}
\beginpicture

\setcoordinatesystem units <9pt,9pt>
\setplotarea x from  -8 to 63, y from -1 to 9

\setlinear \plot 0  0    0 5.8 /
\put {$\circ$} at 0  6
\setlinear \plot 0  6.3    0 8 /

\setlinear \plot 8.66 0   8.66 8 /

\circulararc -120 degrees from 0 0 center at 5.77 0


\circulararc -85 degrees from 0 0 center at 2.885 0
\put {$\circ$} at 2.885  2.885
\circulararc  85 degrees from 5.77 0 center at 2.885 0

\circulararc -60 degrees from 5.77 0 center at 11.55 0

\circulararc -80 degrees from 5.77 0 center at 7.215 0
\put {$\circ$} at 7.21  1.445
\circulararc  80 degrees from 8.66 0 center at 7.215 0

\put { \small    $0$} at   0 -2
\put { \small $1/\lambda$} at   5.77 -2
\put { \small $\lambda/2$} at   8.66 -2

\put {$ e_3(1)$} at   10.5  5



\setlinear \plot 20  0    20 8 /
\put {\small $\vee$} at 20  5
\put {\tiny $ <$} at   27.215 1.445

\setlinear \plot 28.66 0   28.66 8 /

\circulararc -120 degrees from 20 0 center at 25.77 0


\circulararc -85 degrees from 20 0 center at 22.885 0
\put {$\circ$} at 22.885  2.885
\circulararc  85 degrees from 25.77 0 center at 22.885 0

\circulararc -60 degrees from 25.77 0 center at 31.55 0

\circulararc -180 degrees from 25.77 0 center at 27.21 0

\put { \small    $0$} at   20 -2
\put { \small $1/\lambda$} at   25.77 -2
\put { \small $\lambda/2$} at   28.66 -2

\put {$ e_3(1)$} at   30.5  5



\put {\bf Figure 3a} at 4.77  -5
\put {\bf Figure 3b} at  24.77 -5
\endpicture
\end{figure}
\end{center}

\vspace{-.3cm}

\noindent (ii) Let
$ M_Y= \{
   {-\infty} \,\,_{_{_{_{\smile}} }}\  \hspace{-.37cm}_{ _{_{_{_{_{_{_{1}}}}}}} }
   {0} \,\,_{_{_{_{\smile}} }}\  \hspace{-.37cm}_{ _{_{_{_{_{_{_{\circ}}}}}}} }
 {1/\lambda} \,\,_{_{_{_{\smile}} }}\  \hspace{-.37cm}_{ _{_{_{_{_{_{_{1}}}}}}} }
     {\lambda/2 } \,\,_{_{_{_{\smile}} }}\  \hspace{-.52cm}_{ _{_{_{_{_{_{_{e_3(1)}}}}}}} }
  \infty   %
\}$ (see Figure 3b).
 $Y$ is a subgroup of index 3 in $G_6$.
 An easy calculation
      of the side pairings shows that $Y = G_0(2)=\{
       (x_{ij} )\in G_6 : x_{12} \equiv 0\,(mod\,\,2)\}$.
         Note that
         both $M_X$ and $M_Y$ are 3-clusters (with different
          side pairings).


\subsection{Subgroups of index 2 of $G_q$}
 A special polygon of a subgroup of index 2 is either a 2-cluster or a union
  of two special triangles (1-clusters).
 In the case $q$ is odd, it has to be a union of two special triangles. As a
  consequence, its  Hecke-Farey Symbol is given as follows.
  $$
M_1 =
\{-\infty
 _{_{_{\smile}}}\  \hspace{-.36cm} _{_{_{_{_{_{_{_{\bullet}}}}}}}}
 {0}
_{_{_{\smile}} }\  \hspace{-.37cm}_{ _{_{_{_{_{_{_{\bullet}}}}}}}}
 {\infty}
 \}. \eqno(5.3)$$

 \smallskip
 \noindent
 A set of independent generators is $\{ ST^{-1}, T^{-1}S\}$.
 In the case $q$ is even, $G_q$ has three subgroups of index 2, (5.3) and two more
  given as follows.
  $$
M_2 =
\{-\infty
 _{_{_{\smile}}}\  \hspace{-.36cm} _{_{_{_{_{_{_{_{\circ}}}}}}}}
 {0}
_{_{_{\smile}} }\  \hspace{-.37cm}_{ _{_{_{_{_{_{_{\circ}}}}}}}}
{1/\lambda_q} \,\,\,
_{_{_{\smile}} }\  \hspace{-.52cm}_{ _{_{_{_{_{_{_{e_2(1)}}}}}}}}
 {\infty}
 \},
 \,\,\,
 M_3 =
\{-\infty
 _{_{_{\smile}}}\  \hspace{-.36cm} _{_{_{_{_{_{_{_{1}}}}}}}}
 {0}
_{_{_{\smile}} }\  \hspace{-.37cm}_{ _{_{_{_{_{_{_{1}}}}}}}}
{1/\lambda_q} \,\,\,
_{_{_{\smile}} }\  \hspace{-.52cm}_{ _{_{_{_{_{_{_{e_2(1)}}}}}}}}
 {\infty}
 \}.   \eqno(5.4)$$

\smallskip
\noindent To determine $e_2(1)$, we once again follow the notation of
 (i) of subsection 4.3.
 One has $r=2$,  $g=1$,  $x_i= 1/\lambda_q$, and $x_{i+1}= 1/0$.
   It follows from (4.4) that
     $e_2 (1)= R^2$. In particular,
 $$M_3
=\left <
{ R^2=
 \left (
\begin{array}{rc}
0 & 1\\
-1 & \lambda_q \\
\end{array}
\right ) ^2,
 \left (
\begin{array}{cc}
1 & 0\\
 \lambda_q &1\\
\end{array}
\right ) }
\right >.\eqno(5.5)$$
$M_3$ is known as the {\em even subgroup} of $G_q$ ($q$ even).
   Note that subgroups of index 2 of $G_q$ $(q$ even)
    can be realised as Veech groups.

\subsection{Commutator subgroups and principal congruence subgroups of level 2}
 Let $P_0=\{ x_1 , x_2, \cdots, x_{q-1} , x_q  \}$
 be the cusps of the depth 1 $q$-gon
  in the right half plane. The reduced forms of the cusps of $P_0$ can be
  determined by Lemma 3.2. Take note that $x_1=0$ and
   that $x_q =\infty$. Throughout the subsection, the $x_i$'s
   are given as above and that $ x_i <x_{i+1}$.

 \smallskip
 \noindent {\bf Proposition 5.3.}
  {\em A special polygon for $G_q'$, the commutator subgroup of $G_q$, is the
   union of two depth one $q$-gons  forming a polygon with $2(q-1)$ free sides as follows.
  $G_q'\setminus \Bbb H$ is a genus $[(q-1)/2]$
   surface with  one
   cusp
    $($resp. two cusps$)$ if $q$ is is odd $($resp. even$)$.
   $$\{
  -\infty \,\, _{_{_{\smile}}} \hspace{-.22cm}_{ _{_{_{_{_{_{_{1}}}}}}} }
  -x_{q-1} \,\, _{_{_{\smile}}} \hspace{-.22cm}_{ _{_{_{_{_{_{_{2}}}}}}} }
  \cdots
   \,\,\,\, _{_{_{\smile}}} \hspace{-.4cm}_{ _{_{_{_{_{_{_{q-2}}}}}}} }
   -x_2
  \,\, \,\,_{_{_{\smile}}} \hspace{-.4cm}_{ _{_{_{_{_{_{_{q-1}}}}}}} }
  0\,\,
  _{_{_{\smile}}} \hspace{-.22cm}_{ _{_{_{_{_{_{_{1}}}}}}} }
  x_2 \,\, _{_{_{\smile}}} \hspace{-.22cm}_{ _{_{_{_{_{_{_{2}}}}}}} }
  \cdots
   \,\,\,\, _{_{_{\smile}}} \hspace{-.4cm}_{ _{_{_{_{_{_{_{q-2}}}}}}} }
   x_{q-1}
  \,\, \,\,_{_{_{\smile}}} \hspace{-.4cm}_{ _{_{_{_{_{_{_{q-1}}}}}}} }
  \infty
 \}.\eqno(5.6)$$}

\noindent {\em Proof.}
 The special polygon in (5.6) contains $2q$ special triangles. Hence the side
  pairings given by (5.6) generates a subgroup of index $2q$.
As  $[G_q : G_q'] = 2q$, it thus suffices to prove
  that the side pairings given in (5.6) are members in $G_q'$.
   This can be checked easily as
    the side pairing that pairs $(x_i, x_{i+1})$ and
   $(-x_{q-i+1}, -x_{q-i})$ is given by $T^{-1} R^i TR^{-i}$.\qed

\smallskip
\noindent {\bf Proposition 5.4.} {\em Let $q\ge 3$ be a prime.
A special polygon for  $G (2)$, the principal congruence subgroup of level $2$, is the
   union of two depth one $q$-gons  forming a polygon with $2(q-1)$ free sides as follows.
   $G(2)\setminus \Bbb H$ is a genus zero surface with $q$ cusps.
    %
  $$\{
  -\infty \,\, _{_{_{\smile}}} \hspace{-.22cm}_{ _{_{_{_{_{_{_{1}}}}}}} }
  -x_{q-1} \,\, _{_{_{\smile}}} \hspace{-.22cm}_{ _{_{_{_{_{_{_{2}}}}}}} }
  \cdots
   \,\,\,\, _{_{_{\smile}}} \hspace{-.4cm}_{ _{_{_{_{_{_{_{q-2}}}}}}} }
   -x_2
  \,\, \,\,_{_{_{\smile}}} \hspace{-.4cm}_{ _{_{_{_{_{_{_{q-1}}}}}}} }
  0\,\,\,\,
  _{_{_{\smile}}} \hspace{-.4cm}_{ _{_{_{_{_{_{_{q-1}}}}}}} }
  x_2\,\, _{_{_{\smile}}} \hspace{-.4cm}_{ _{_{_{_{_{_{_{q-2}}}}}}} }
  \cdots
   \,\,\, _{_{_{\smile}}} \hspace{-.22cm}_{ _{_{_{_{_{_{_{2}}}}}}} }
   \,\,\,x_{q-1}
  \,\, \,_{_{_{\smile}}} \hspace{-.22cm}_{ _{_{_{_{_{_{_{1}}}}}}} }
  \infty
 \}.\eqno(5.7)$$}

\noindent {\em Proof.}   Since $q$ is a prime,  $[G_q : G(2)] = 2q$.
By (4.7),  the side pairing that pairs $(x_i, x_{i+1})$ and $(-x_{i+1}, -x_i)$
  is  in $G(2)$ for every $i$.
   This completes the proof our assertion.   \qed

\smallskip
 In the case $q\ge 3$ is a prime,
similar to Proposition 5.4, one can show that $ G_0(2)=G_1(2)$ has index $q$ and admits
 the following special polygon $(P, I_X)$, where $P$ is a $q$-gon.
$$\{
  -\infty\,
 _{_{_{\smile}}} \hspace{-.22cm}_{ _{_{_{_{_{_{_{1}}}}}}} }
  x_1\,
  _{_{_{\smile}}} \hspace{-.22cm}_{ _{_{_{_{_{_{_{2}}}}}}} }
  x_2\,\,
  \cdots\,\,
  x_{(q-1)/2}
  \,\, _{_{_{\smile}}} \hspace{-.22cm}_{ _{_{_{_{_{_{_{\circ}}}}}}} }
  x_{(q+1)/2}\,\,\,\,
  %
  \cdots
   \,\,\, _{_{_{\smile}}} \hspace{-.22cm}_{ _{_{_{_{_{_{_{2}}}}}}} }
   \,\,\,x_{q-1}
  \,\, \,_{_{_{\smile}}} \hspace{-.22cm}_{ _{_{_{_{_{_{_{1}}}}}}} }
  \infty
 \}.\eqno(5.8)$$

 $G_0(2) \setminus \Bbb H$
 is a genus zero surface with one elliptic element of order 2 and $(q+1)/2$ cusps.
 Note that the above is not true if $q$ is not a prime.
 For instance, the group in (ii) of subsection 5.1 gives $G_0(2) $ of $G_6$.
  Its special polygon is a 3-cluster, not a 6-gon.

\subsection{Power subgroups}
 For each $n\in \Bbb N$, the power subgroup $G_q^n$
 is the characteristic  subgroup of $G_q$ generated by $\{x^n\,:\, x\in G_q\}$.
Study of power subgroups has a long history back to Newman [N].

\smallskip
\noindent {\bf Proposition 5.5.} {\em Let $q\ge 3$ be an odd integer.
 Then $G_q^2$ is a normal subgroup of index $2$. A
  special polygon of $G_q^2$ is a union of two special triangles
   given as in $(5.3)$. In the case $q$ is even,
    a Hecke-Farey symbol of $G_q^2 $ is given as in $M_3$ of $(5.4)$.}

\smallskip
\noindent {\em Proof.} We shall first assume that $q$ is odd.
 It is clear that $A = T^{-1} S$ (given as in Lemma 3.2) and $R = ST^{-1}$,
  each has order $q$, are elements of $G_q^2$. Hence the special polygon of  (5.3)
   is also a special polygon of $G_q^2$.

    In the case $q$ is even,  we consider the homomorphism $\phi \, :\, G_q\to
     \Bbb Z_2 =\left < a\right >$ defined by $\phi (S) = \phi (R) = a$.
      It is clear that $G_q^2$ is the kernel of $\phi$.
       Since $G_q$ has exactly three subgroups of index 2 and
        the conjugates of $S$  and $R$ cannot be members of
         a set of independent generators of $G_q^2$,
          we conclude that
         $M_1$ and $M_2$ (see (5.3) and (5.4)) are not   special polygons of $G_q^2$. As a consequence,
           a  special polygon of
           $G_q^2$ is given as in $M_3$ of (5.4)\qed

\smallskip
\noindent {\bf Proposition 5.6.} {\em
Let $\{  x_1, x_2, \cdots, x_{q-1}, \infty\}$ be given as in subsection $5.3$
 and let $r>1$ be an odd divisor of $q$. Then
$[G_q : G_q^r] = r$ and
a Hecke-Farey symbol of $G_q^r$ is given
by}
   $$\{
  -\infty \,\, _{_{_{\smile}}} \hspace{-.22cm}_{ _{_{_{_{_{_{_{\circ}}}}}}} }
   x_1\,\,
  _{_{_{\smile}}} \hspace{-.22cm}_{ _{_{_{_{_{_{_{\circ}}}}}}} }
  x_2 \,\, _{_{_{\smile}}} \hspace{-.22cm}_{ _{_{_{_{_{_{_{\circ}}}}}}} }
  \cdots
   \,\,\,\, _{_{_{\smile}}} \hspace{-.22cm}_{ _{_{_{_{_{_{_{\circ}}}}}}} }
   x_{r}
  \,\, \,\,_{_{_{\smile}}} \hspace{-.37cm}_{ _{_{_{_{_{_{_{e_r(1)}}}}}}} }
  \infty
 \}.\eqno(5.8)$$
\noindent {\em Proof.}
Define $\phi \,:\, G_q \to \Bbb Z_{r}= \left < a\right >$ by $\phi(S)= 1$, $\phi(R) =
 a$. It is clear that $\phi$ is a homomorphism and that $G_q^r \subseteq \, ker\phi$.
Hence $[G_q : G_q^r]
\ge r$.
 Let $G$  be the group generated by the side pairings of (5.8).
The  side pairings of (5.8) are conjugates of $S=S^r$ and $e_r(1) = R^r$.
 Hence
 $G \subseteq G_q^r$.
  Since $[G_q : G]= r$, one has $[G_q : G_q^r] \le r$.
   Hence $G_q^r = G$.\qed


\subsection{Non-free normal subgroups} Let $X$ be a proper normal subgroup of $G_q$ that
 is not a free group. It follows that either $R\in X$ or $S\in X$.
  Suppose that $R\in X$. Since $X$ contains all the conjugates of $R$,
   both $R= ST^{-1}$ and $T^{-1}S$ are members of $X$. Hence $X$ has index 2
    and a special polygon of $X$ is given as in (5.3).
     Note that $X$ is a free product of two copies of $\Bbb Z_q$.
      In the case $S\in X$, $G_q/X$ can be generated by $RX$.
       Since $R$ has order $q$, $G_q/X$ is a cyclic group of order $r$, where $r|q$.

    \section{ Hurwitz-Nielsen realisation problem }

\subsection {Geometric invariants} An immediate application of the study of the special polygons is that
 the geometric invariants of $X\setminus \Bbb H$ can be determined easily.
 Let $M_X= (P,I_X)$ be a special polygon  associated with $X\subseteq G_q$.
 \begin{enumerate}

 \item[(i)] $[G_q : X]=$  the number of special
  triangles in  $M_X= (P, I_X) $.

  \item[(ii)] The subgroup $X$ has
   $\tau _2$
   (the number of the circles $\circ$ in $M_X$) inequivalent classes of elliptic
    elements of order 2  that are conjugates of $S$.
     \item[(iii)] The subgroup $X$ has
   $v_q$
   (the number of the bullets $\bullet$ in $M_X$) inequivalent classes of elliptic
    elements of order $q$ that are conjugates of $ R= ST^{-1}$.
     \item[(iv)] Let $r$ ($1< r<q$) be a divisor of $q$. $X$ has
   $v_r$
   (the number $e_{q/r}(g)$'s in $M_X$) inequivalent classes of elliptic
    elements of order $r$ that are conjugates of $ R^{q/r}$.
     \item[(v)]
      Suppose that the
       cusps of $P$
     is  partitioned into $v_{\infty}$  classes under the action of
     $I_X$.
      Then $v_{\infty} = $ the number of cusps of $X\setminus \Bbb H = $
         the number
       of cycles of $f(T)$,  where
      $f(T)$ is the permutation representation of $T$ on
      the set of cosets  $ G_q/X$.

 \item[(vi)] Let $\Delta =\{ r\,:\, r\mbox{ is a divisor of } q, \,2\le r\le q\}
  = \{r_1, r_2, \cdots, r_k=q\}$. The genus $g$ of $X\setminus \Bbb H$ is given by the following
   Riemann Hurwitz formula.
$$ 2g-2 + \tau_2/2 + \sum_{i=1}^k v_{r_i}(1-1/r_i) + v_{\infty} = [G_q : X] (1/2-1/q)
.\eqno(6.1)$$


\end{enumerate}

\noindent The terms in  $\{g, \tau_2, v_{r_1}, v_{r_2},\cdots, v_{r_k}, v_{\infty}, [G_q : X]\}$ are called the {\em geometric invariants}
 of $X$, where the $v_i$'s are given as in (ii)-(iv).

 \smallskip
 \noindent {\em Proof.} (i)  is clear.
  (ii)-(iv) follows from the observation that
  each conjugacy class of elliptic elements of $X$ must have
   exactly one representative in $I_X$ as $I_X$ is a set of independent
    generators. (v) follows from the fact that the number of cusps is
     just the number of double cosets $\left <T\right > \setminus G_q\,/\,X$.
 \qed

\smallskip
\noindent {\bf Lemma 6.1.} {\em
Suppose that a special polygon $M_X=(P,I_X)$ of $X$ consists of $n_0$
 $q$-gons
  and $v_{r_i}$ $q/r_i$-clusters, $1\le i \le k$.
   Then the  number of generators of infinite order in $I_X$ is
$$
 f = \frac{n_0(q-2) + \sum_{i=1} ^ k v_{r_i}(q/r_i-2) + 2- \tau_2}{2}.\eqno(6.2)$$}
\noindent {\em Proof.}
          Since the boundary of a $q/r$-cluster has $q/r$ even lines,
       it follows  that  the boundary of $P$ has
     $n_0(q-2) + \sum_{i=1}^k v_{r_i} (q/r_i-2)  +2$ even lines.
     Among these even lines, $\tau_2$ of them are self paired by elements of $I_X$.
As a consequence,
the remaining even lines are free sides and they  are paired
          by side pairings of infinite order. In particular,
           the number of such side pairings is
           $ f = (n_0(q-2)  +\sum_{i=1 }^kv_{r_i}(q/r_i-2) +2- \tau_2 )/2.$\qed

\subsection{Construction of convex hyperbolic polygons}
 Let $P_0$ be either the depth one $q$-gon in the right  half plane or the
  $r$-cluster  $\Phi_r$
 where $r\ge 2$ is a divisor of $q$ (see (4.1a) for $\Phi_r$). The Hecke-Farey sequence associated with $P_0$ is given
  by
  $$F_0=\{-\infty, 0, a_1, a_2, \cdots, a_{k-1}, a_k = \infty\}. \eqno(6.3)$$
Since $r\ge 2$, $\{0, a_1\}$ is an ordinary interval of $F_0$.
In the case $P_0$ is $\Phi_r$,
 $\{a_{k-1}, a_k\}$ is an $r$-interval of $F_0$ (see (d) of subsection 4.3 for the terms
  ordinary and $r$-intervals).
 Note that   $P_0$ is convex.
  Since  $(0, a_1)$
   is an even line of $P_0$ (see (d) of subsection 4.3), one may attach $R_0$
    to $P_0$ along the even line $(0, a_1)$
    to get a new polygon $P_1$, where $R_0$ is
     either a $q$-gon or an $s$-cluster $(s \ge 2$ and $s|q)$.
     The Hecke-Farey sequence associated with $P_1$ takes
     the form
     $$F_1 = \{-\infty, 0, b_1, b_2, \cdots, b_{t-1}, a_1, a_2, \cdots, a_k\}.\eqno(6.4)$$

\noindent Since $s\ge 2$,  $\{0, b_1\}$ is an ordinary interval of $P_1$ and $\{b_{t-1}, a_1\}$
is an $s$-interval if $R_0$ is an $s$-cluster. Note that  $P_1$ is convex.
 Since $(0, b_1)$ is an even line of $P_1$,
 one may attach $R_1$
    to $P_1$ along the even line $(0, b_1)$
    to get a new polygon $P_2$, where $R_1$ is
     either a $q$-gon or an $u$-cluster $(u \ge 2$ and $u|q)$.
 Apply this procedure repeatedly, one admits a convex hyperbolic polygon
 $P_n$ consists of $n_0$ $q$-gons and  $v_{r_i}$ $q/r_i$-clusters for some $n_0$
  and  $v_{r_i}$, where
  $r_i|q$ and $q/r_i \ge 2$
(the key of our construction is that $P_i$  always has an even line $(0, x_i)$
 and that $R_i$ is always attached to $P_i$ along $(0, x_i)$).
    Let $\Bbb F$ be the Hecke-Farey sequence associated with $P_n$.
     $\Bbb F$ has $v_{r_i}$  $q/r_i$-intervals, where $q/r_i \ge 2$.
      The intervals of
 $\Bbb F$ are divided into two classes (i)  $q/r_i$-intervals, where $q/r_i \ge 2$,
  and (ii) ordinary intervals.
   To count the number of ordinary intervals, we note that
\begin{enumerate}
\item[(i)] the boundary of a $q$-gon has $q$ even line,
 \item[(ii)]  the boundary of a $q/r_i$-cluster has $q/r_i$ even lines.
\end{enumerate}

\noindent It follows that the boundary of $P_n$ has $n_0(q-2)
 + \sum v_{r_i}(q/r_i-2) +2 $ even lines. As the even lines of $P_n$ are associated
  with the ordinary intervals of $\Bbb F$ (see (d) of subsection 4.3),
  the  number of ordinary intervals of $\Bbb F$ is  given
   by
 $$n_0(q-2)
 + \sum v_{r_i}(q/r_i-2) +2, \,\, q/r_i \ge 2. \eqno(6.5)$$

 \subsection{Millington's Theroem}
  Kulkarni [K2, K3] gives two  proofs of Millingtion's Theorem [Mi], one by
   {\em Diagrams} (see Section 4 of [K2]) and one by
   {\em Farey symbols} (see Section 7.6 of [K3]). We extend Millington's result to $G_q$
     by studying  Hecke-Farey symbols. Our proof is  a simple generalisation
     of Kulkarni's proof of  Millingtion's Theorem.

 \smallskip
\noindent {\bf Theorem 6.2.} {\em
Let  $g\ge 0, \, \tau_2\ge 0,\, v_{r_i}\ge 0$, $d\ge 1,\, v_{\infty} \ge 1$
be integers and let $r_i \in \Delta =\{ r\,:\, r\mbox{ is a divisor of } q, \,2\le r\le q\}
  = \{r_1, r_2, \cdots, r_k=q\}$.
 Then $G_q$ has a subgroup $X $ of index $d$, genus $g$ with
  $v_{\infty} $ cusps, and $\tau_2$ $(resp.\,\, v_{r_i})$ conjugacy classes of
 elements that are conjugates of $S$ $(resp.\, \,R^{q/r_i})$
 if and only if $ m_0= 4g-4 +\tau_2 + 2v_{\infty} + \sum _{i=1}^k v_{r_i}(2-q/r_i)\ge 0$ is a multiple of $(q-2)$  and
  $$ 2g-2 + \tau_2/2 +\sum_{i=1}^k  v_{r_i}(1-1/r_i)  + v_{\infty} = d(1/2-1/q)
.\eqno(6.6)$$
Note that if $q$ is odd,   then $m_0/(q-2) \in \Bbb Z$ is a consequence of $(6.6)$. }

\smallskip

\noindent {\em Proof.}
Suppose that $m_0\ge 0$ is a multiple of $q-2$ and that
 $g\ge 0, \, \tau_2\ge 0,\, v_{r_i}\ge 0$, $d\ge 1,\, v_{\infty} \ge 1$ satisfy the Riemann-Hurwitz formula (6.6).
Let $r_i\in \Delta$. It follows that $q/r_i \ge 2$ for $1\le i\le k-1$. Let $n_0=m_0/(q-2)$.

\smallskip
\noindent {\bf Case 1.  $n_0=0$ and $v_{r_i }=0$ for all $i\le k-1$.}
 A simple calculation shows that either  (i) $d=1, \tau_2=v_q =1$, $v_{\infty} = 1$, $g=0$, $ v_{r_i}=0$ for all $i\le k-1$
   or (ii) $d=2, v_q =2$, $v_{\infty} = 1$, $g=0$,  $ v_{r_i}=0$ for all $i\le k-1$.
     In case (i), $X = G_q$.  In case (ii), $X$ is given as in (5.3).

\smallskip
\noindent {\bf Case 2.  $n_0>0$ or  $v_{r_i }>0$ for some $i\le k-1$.}
By (6.5), we have a polygon $P_n$ and a Hecke-Farey sequence $\Bbb F$  with
 $v_{r_i}$ $q/r_i$-intervals $(i= 1,2,\cdots,\, k-1)$  and
 $$n_0(q-2)
 + \sum_{i=1}^{k-1} v_{r_i}(q/r_i-2) +2
  = 4g-2+\tau_2 +v_q + 2v_{\infty}\ge 2 \eqno(6.7)$$

\noindent ordinary intervals.
 Since $n_0>0$ or  $v_{r_i} >  0$ for some  $i \le k-1$,
  $\Bbb F$ has at least three intervals.
  We make $\Bbb F$ into a Hecke-Farey symbol by declaring
  the first  $\tau_2$  ordinary intervals
    even intervals (see Definition 4.2 for the terms even and odd intervals)
  and
  the next
   $v_q$ ordinary intervals  odd intervals (this is equivalent to adjoin $v_q$ special
    triangles to $P_n$). The next $2(v_{\infty}-1)$  ordinary
    intervals    are declared to be free
     intervals and are divided into $(v_{\infty}-1)$ consecutive pairs which are
      paired. The remaining $4g$  ordinary intervals  are free intervals and are paired
       in the usual $xyx^{-1}y^{-1}$ fashion (see Discussion 6.3).
  Let $X$ be the group generated by the above side pairings and let  $I_X$ be the set of
   the side pairings.
   \begin{enumerate}
   \item [(i)] $I_X$ possesses  $\tau_2$  elements that are conjugates of $S$ and
    $v_{r_i}$ elements that are  conjugates of $R^{q/r_i}$ for $1\le i\le k$.

   \item[(ii)] An easy study of the side pairings
    shows that $X$ has $v_{\infty}$ cusps. The index of $X$ is given
    by the number of special triangles of the special polygon $M_X= (P, I_X)$,
     where $P$ is the union of $P_n$ and $v_q$ special triangles.
    There are $n_0=m_0/(q-2)$ $q$-gons,
     $v_{r_i}$ $q/r_i$-clusters $(1\le i\le k-1)$ and $v_q$ special triangles.
     Hence
     the index is
    $$n_0 q + \sum_{i=1}^{k-1} v_{r_i} (q/r_i) +v_q  =d.\eqno(6.8)$$
  \end{enumerate}

 \noindent  The genus of $X$ is $g$ (see (6.1)).
    It follows that  $X$ is the required subgroup (Theorem 5.1).

\smallskip
 Conversely,  let $X$ be given as in our theorem,   the geometric
 invariants of $X$ satisfy (6.6). Let $M_X=(P, I_X)$ be a special polygon of $X$.
  Then $P$ has $v_{r_i}$ $q/r_i$-clusters $(1\le i \le k$). Suppose that $P$ has $n_0$ $q$-gons.
   By (i) of subsection 6.1,
    $d = n_0q +\sum_{i=1}^k v_{r_i} q/r_i$. By (6.6), one has
    $n_0 =
   (4g-4 +\tau_2 + 2v_{\infty} + \sum _{i=1}^k v_{r_i}(2-q/r_i))/(q-2)$.
    In particular,
 $4g-4 +\tau_2 + 2v_{\infty} + \sum _{i=1}^k v_{r_i}(2-q/r_i)\ge 0$ is a multiple of
  $q-2$.\qed

\smallskip
\noindent {\bf Discussion 6.3.}
(i) A Hecke-Farey symbol must have at least two side pairings (subsection 4.3).
(ii) Since $\Bbb F $ has at least three intervals,  our declaration gives
 at least two side pairings.

\subsection{Kurosh's Theorem} 
Let $Y$ be a free product of $F_f$ (a free group of rank $f$), $\pi_2$
    copies of $\Bbb Z_2$, and $v_{r_i} $ copies of $\Bbb Z_{r_i}$, where
    $r_i \in \Delta_0 =\{ r\,:\, r\mbox{ is a divisor of } q, \, 3\le r\le q\}
  = \{r_1, r_2, \cdots, r_s=q\}$.
    We say $Y$ is {\em realisable } in $G_q$
      if $Y$ is isomorphic to a subgroup $X$ of finite index of  $G_q$.
       It is clear that $\Bbb Z_2 *\Bbb Z_q$ and $\Bbb Z_q *\Bbb Z_q$
        are realisable (see (5.3) for $\Bbb Z_q *\Bbb Z_q$).

  \smallskip
  \noindent {\bf Proposition 6.4.} {\em
   Let $q\ge 3$ be an odd integer and let $Y$ be
   given as  above. Suppose that $Y$ is not  $\Bbb Z_2 *\Bbb Z_q$ or
   $\Bbb Z_q *\Bbb Z_q$. Then $Y$ is isomorphic to a subgroup $X$ of  $G_q$ if and only if
    $2f+\pi_2 +v_q -2 > 0$ and
    $m_0= 2f+ \pi_2  + \sum_{i=1}^s v_{r_i} (2-q/r_i)-2 \ge 0 $ is a multiple of $q-2$.
 The index of $X$ is
  $d= m_0q/(q-2)+\sum  _{i=1}^ s v_{r_i} q/r_i$.

}

 \smallskip

 \noindent
 {\em Proof.}
  Suppose that
 $Y$ is isomorphic to a subgroup $X$ of $G_q$.
  Let $M_X = (P, I_X)$ be  a special polygon of $X$. Then
  $I_X$ has $f$ elements of infinite order, $\pi_2$ elements that are conjugagtes of $S$,
   and
   $v_{r_i}$ elements that are conjugates of $R^{q/r_i}$.
   By Lemma 6.1, $f$, $v_{r_i}$, $n_0\ge 0$,  and $\pi_2$ satisfy the identity given in (6.2). This implies that
    $m_0 = 2f+ \pi_2  + \sum_{i=1}^s v_{r_i} (2-q/r_i)-2 = n_0(q-2)\ge 0$.
     In particular, $m_0\ge 0$ is a multiple of $q-2$.

\smallskip
    Note that $n_0 = m_0/(q-2)$ is the number of $q$-gons of $P$.
      Since $P$ has $n_0$ $q$-gons and $v_{r_i} $ $q/r_i$-clusters,
       the  index of $X$ is
    $d= m_0q/(q-2)+\sum _{i=1}^s  v_{r_i} q/r_i$ (see (i) of subsection 6.1).
      Finally, since  $m_0 = 2f+ \pi_2 +v_q -2  + \sum_{i=1}^{s-1} v_{r_i} (2-q/r_i) = n_0(q-2)\ge 0$
   and  $ (2-q/r_i) <0$ for $i\le s-1$, one has
   $2f+ \pi_2 +v_q -2 \ge 0$.
Suppose that
$2f+\pi_2 +v_q -2  = 0$.
It follows easily that $v_{r_i} = 0$ for  $i\le s-1$.
A simple study shows that
$Y$ is isomorphic to a subgroup $X$ of $G_q$
 only if
  $Y\cong \Bbb Z_2 *\Bbb Z_q$ or $\Bbb Z_q *\Bbb Z_q$.
   A contradiction. Hence $2f+\pi_2 +v_q -2  > 0$.

\smallskip
Conversely, suppose that $2f+\pi_2 +v_q-2 >0$ and that $n_0 = m_0/(q-2)\in \Bbb N \cup\{0\}$.
 If  $n_0 = 0$ and $v_{r_i}=0$ for  $i\le s-1$,
        then $m_0 = 2f+ \pi_2 +v_q -2 =0$. A contradiction. Hence
         either $n_0>0$ or $v_{r_i} >0$ for some $i\le s-1$.
         By (6.5), we have a polygon $P_n$ and a Hecke-Farey sequence $\Bbb F$ with
      $v_{r_i}$ $q/r_i$-intervals $(q/r_i\ge 2$ for $i\le s-1$) and
 $$n_0(q-2)
 + \sum_{i=1}^{s-1} v_{r_i}(q/r_i-2) +2
  = 2f+ \pi_2 +v_q \ge 2  \eqno(6.9)$$

      \noindent ordinary intervals.
      We make $\Bbb F$ into a Hecke-Farey symbol by declaring the first
   $\pi_2$ ordinary intervals  even  intervals and the next $v_q$ ordinary
    intervals odd intervals (this is equivalent to adjoin $v_q$ special
    triangles to $P_n$).
   The last  $2f$  ordinary
    intervals    are declared to be free
     intervals and are divided into $f$ consecutive pairs which are
      paired.
       Let $X$ be the subgroup of $G_q$ generated by the
        above side pairings and let  $I_X$ be the set of the side pairings. Then

     \begin{enumerate}
   \item [(A)] $I_X$ possesses  $\tau_2$  elements that are conjugates of $S$
    (of order 2) and
    $v_{r_i}$ elements that are  conjugates of $R^{q/r_i}$ (of order $r_i$)
    for $1\le i\le s$.

\item[(B)] $I_X$ has $f$ elements of infinite order.
\end{enumerate}

\noindent  By Theorem 5.1,  $X$ is  isomorphic to $Y$. This completes the proof of the
 proposition. \qed

\smallskip
  In the
 case $q$ is even,  subgroups of $G_q$  may possess two types of elliptic elements
  of order 2, namely, conjugates of $S$ and   $R^{q/2}$.
    Note that  $2 \notin \Delta_0$.  For our
         convenience, we set $v_2=0$.

  \smallskip
  \noindent {\bf Proposition 6.5.} {\em
   Let $q\ge 4$ be an even integer and let $Y$ be
   given as  above. Suppose that $Y$ is not  $\Bbb Z_2 *\Bbb Z_q$ or
   $\Bbb Z_q *\Bbb Z_q$.
    Then $Y$ is isomorphic to a subgroup $X$ of  $G_q$ if and only if

   \begin{enumerate}
   \item[(i)]
    $m_0= 2f+ (\pi_2 -t) + t(2-q/2) +\sum_{i=1}^s v_{r_i} (2-q/r_i)-2 \ge 0 $ is a multiple of $q-2$
     for some $t\ge 0$, where $(\pi_2- t)$ is nonnegative, and
   \item[(ii)]   $  2f+(\pi_2 -t)+v_q -2 \ge 0$ and
      $ (2f+(\pi_2 -t)+v_q -2, t,v_{q/2})\ne (0,0,0 )$.

     \end{enumerate}
Note that $v_{2}=0$.
The index of $X$ is
  $d= m_0q/(q-2)+ tq/2 + \sum _{i=1}^s v_{r_i} q/r_i.$

}

\smallskip
 \noindent
 {\em Proof.}
  Suppose that
 $Y\cong X \subseteq G_q$.
  Let $M_X = (P, I_X)$ be  a special polygon of $X$. Then
  $I_X$ has $f$ elements of infinite order, $(\pi_2 -t)$ elements
  that are conjugates of $S$, $t$ elements that are conjugates of $R^{q/2}$, and
   $v_{r_i}$ elements that are conjugates of $R^{q/r_i}$.
    Note that $r_i \ne 2$ ($r_i\in \Delta_0$).
    The  assertion now follows by applying the proof of Proposition 6.4. \qed

\section{Permutation representation of $G_q$ on $ G_q/X $}

The main purpose of this section is to give an easy and systematic
  method that determines the permutation representation
 of $G_q$ on $  G_q/X$.

\subsection{A commutative diagram} Let $\Phi$ be given as in  subsection 4.1.
 $\Phi$ is a fundamental domain of $G_q$.
  For each coset $Xg\in G_q/X$, there exists a unique $x_g\in X$
  such that $x_g g\Phi \in M_X = (P, I_X)$.
  Denoted by $\Omega_X$ the set of all such special triangles $x_g g\Phi$.
  It follows that
     $\cup_{x_g g\Phi \in \Omega_X} x_g g\Phi   = M_X.$
 As a consequence, there is a one to one correspondence between $  G_q/X$ and
 $ \Omega_X$ defined by  $$\tau(Xg) = x_gg\Phi \in \Omega_X.\eqno(7.1)$$

\smallskip

An element $g$ of $G_q$
     acts on $ G_q/X$ by $(X h, g) \to   Xhg$.
We shall now study the action of $X$ on $\Omega _X$ as follows.
     Let $g\in G_q$.
   For each $g_i\Phi\in \Omega _X$,
   there exists a unique pair
    $(g_j\Phi ,  x_{ij})\in \Omega_X \times X$
      such that  $x_{ij}g_ig\Phi =g_j\Phi \in \Omega _X $. The action of $g$ on $\Omega_X$ is defined by
      $$(g_i\Phi, g) \to   x_{ij}g_i g\Phi.\eqno(7.2)$$

  \smallskip
  \noindent {\bf Lemma 7.1.} {\em
  The action of $G_q$ on $G_q/X$ is isomorphic
   to the action of $G_q$ on $\Omega_X$.}

   \smallskip
   \noindent {\em Proof.}
 We consider the following  diagram,
 where the horizontal arrows represent the actions of $X$ on $G_q/X$ and $\Omega_X$.
    $$
    \begin{array}{ccl }
       G_q/X \times  G_q  &  \longrightarrow &  G_q/X\\
(\tau, id) \phantom {\Bigg |}\downarrow\,\,\,\,\,\,\,\,\, && \,\, \downarrow  \tau\\
   \Omega_X\times  G_q  & \longrightarrow  &  \Omega _X\\
     \end{array} \eqno(7.3)$$

\smallskip
\noindent
  Let $G_q /X = \cup Xg_i$ be chosen such that $\cup g_i\Phi = M_X$.
  This implies that $\tau (Xg_i) = g_i\Phi$. As a consequence, (7.3)
   is a commutative diagram and the action of $G_q$ on $G_q/X$
    is isomorphic to the
     action of $G_q$ on $\Omega_X$.\qed

  \subsection{The permutation representation of $G_q$ on $\Omega_X$}

   The main purpose of this subsection is to give the permutation representations
    of $S$ and $R$ on $\Omega _X$.
    \smallskip

  \noindent {\bf Lemma 7.2.} {\em Let $\Omega_X = \{ g_i\Phi\,:\, \cup g_i \Phi= M_X\}$.
 The permutation representation of $S$ on $\Omega _X$
is given by
$$f(S) = \prod (g_i\Phi, g_j\Phi),\eqno(7.4)$$

\noindent
where $g_i\Phi\in \Omega _X$ and $g_j\Phi\in \Omega_X$ share the same even line of $P$.
$(g_i\Phi, g_j\Phi)$ is a one cycle if and only if
$g_i S g_i^{-1}\in X$   if and only if $g_i(0,\infty)$  is paired
 with itself by $g_i S g_i^{-1}\in X$.
}

\smallskip
\noindent {\em Proof.}
Let $g_i\Phi \in \Omega_X$. There exists a unique pair $(g_j\Phi, x_{ij})\in \Omega_X \times
 X$ such that
 $x_{ij} g_i S\Phi  = g_j\Phi \in \Omega_X$ (see (7.2)). By Lemma 7.1,
the permutation representation of $S$ is given
by $\prod (\tau(Xg_i), \tau(Xg_i S))=
\prod (g_i\Phi , x_{ij} g_iS\Phi)$
$=\prod (g_i\Phi,  g_j\Phi)$.
The even line of $g_j\Phi \in M_X$ is
$  g_j(0,\infty) = x_{ij} g_i S (0,\infty) =x_{ij} g_i(0,\infty) \in M_X$.
Since $ g_i(0,\infty)\in M_X$, $x_{ij}\in X$ and every even line has a unique $X$-image
 in $M_X$, we conclude
  that $g_j(0,\infty) =x_{ij}g_i(0,\infty)  = g_i(0,\infty)$.
   Hence $g_j\Phi$ and $g_i\Phi$ share the same even line.

\smallskip
 $(g_i\Phi, g_j\Phi)$ is a one cycle if and only if
$(Xg_i, Xg_iS)$ is a one cycle (see Lemma 7.1 and  (7.1)) if and only if $g_iSg_i^{-1}\in X$
 if and only if  $g_i(0,\infty)$ is paired with itself by
 $g_iSg_i^{-1}\in X$.  \qed

  \smallskip

\noindent {\bf Discussion  7.3.}
Let
$ E= E_1\cup E_2$ be the set of even lines of $M_X= (P, I_X)$,
where    $E_1
=\{L_1,L_2,\cdots L_m\}
$ is the set of even lines that are not paired with itself
 by  elements of $ I_X$ of order 2, and
 $E_2
 =\{L_{m+1}, \cdots, L_n\} $ is the set of even lines that
 are paired with itself by
   elements of $I_X$ of order 2.
  Consequently, an element $L_s$ of $E_2$ belongs to exactly one
special triangle of  $\Omega_X$.
  Lemma 7.2 suggests a simple way to list the special triangles in $\Omega_X$. Namely,
 the special triangles share the even line $L_r\in E_1$ are labeled as $r$ and $\overline
  r$ and  the special triangle contains  $L_s\in E_2$  is labeled as $s$.

  \smallskip

  \noindent {\bf Example 7.4.}
 Let
$$ M_X = \{  {-\infty }_{_{_{\smile}} }\  \hspace{-.37cm}_{ _{_{_{_{_{_{_{\bullet}}}}}}} }
 0
_{_{_{\smile}} }\  \hspace{-.37cm}_{ _{_{_{_{_{_{_{\circ}}}}}}}}
1\,\,
 _{_{_{\smile}}}\  \hspace{-.37cm} _{_{_{_{_{_{_{_{1}}}}}}}}
 2\,\,
_{_{_{\smile}} }\  \hspace{-.37cm}_{ _{_{_{_{_{_{_{1}}}}}}}}
3
 _{_{_{\smile}}}\  \hspace{-.36cm} _{_{_{_{_{_{_{_{\bullet}}}}}}}}
 \infty
  \}\eqno(7.5)$$
be the Hecke-Farey symbol of a subgroup $X$ of index 11 of $G_3 = PSL(2, \Bbb Z)$
 (see Figure 4).
 $M_X$ has
 6 even lines, five of them are shared by two special triangles of $\Omega_X$,
  they are
 $$L_1 =(0, \infty), L_2= (1, \infty), L_3 =(2, \infty), L_4= (3, \infty),
 L_5 = (1,2).\eqno (7.6)$$
 Note that  $(1,2)$ and $(2,3)$ give  the same line as they are paired by the side pairing labeled by the
  natural number $1$ (see (7.5)).
 The even line $ L_6=(0,1)$ is paired to itself by an element of order 2 (see (7.5)).
  As a consequence, the special triangles of $\Omega _X$ are labeled as in Figure 4.
 By Lemma 7.2, $f(S) $ is given by
 $$ f(S)= (1, \bar 1)(2,\bar 2)(3, \bar 3)(4, \bar 4)(5, \overline 5)(6).\eqno(7.7)$$

\allowdisplaybreaks

\begin{center}
\begin{figure}
\beginpicture

\setcoordinatesystem units <4pt,4pt>

\setplotarea x from   -26.5 to 83, y from -1 to 22

\circulararc -180 degrees from 0 0 center at 7.5 0
\circulararc -180 degrees from 15 0 center at 22.5 0
\circulararc -180 degrees from 30 0 center at 37.5 0

\circulararc -60 degrees from 45 0 center at  60 0
\circulararc  60 degrees from 0 0 center at  -15 0


\setlinear \plot 0 0 45 0 /

\setlinear \plot  0   0 0 20 /
\setlinear \plot 15  0 15 20 /
\setlinear \plot 30  0 30 20 /
\setlinear \plot 45  0 45 20 /

\setlinear \plot -7.5  13 -7.5 20 /
\setlinear \plot 52.5  13  52.5 20 /

\put {$2$} at   17 14
\put {$3$} at  32  14

\put {$1$} at   2  14
\put {$\overline 4$} at  43  14

\put {$\overline 2$} at   13 14
\put {$\overline 3$} at  28  14

\put {$\overline 1$} at   -2  14
\put {$ 4$} at  47  14

\put {$6$} at  7.3 9.5

\put {$\overline 5$} at  37.5 9.5
\put {$5$} at  22.5  9.5


\put {$\scriptsize{0}$} at 0 -3
\put {$\scriptsize{1}$} at 15 -3
\put {$\scriptsize{2}$} at 30 -3

\put {$\scriptsize{3}$} at 45 -3

\put {\bf Figure 4} at 22.5 -7
\endpicture
\end{figure}
\end{center}

\vspace{-.5cm}

 \smallskip
 The permutation representation of $R$ on $G_q/X$ that contains  $Xg_i\in \Omega_X$ is given by
  the cycle
 $(Xg_i, Xg_iR, \cdots , Xg_iR^{r-1})$,
  where $r$ is the smallest positive integer such that
   $Xg_i = Xg_i R^r$.
  By Lemma 7.1, the permutation representation of $R$ on $\Omega _X$
   that contains $g_i\Phi$
   is given by  the following $r$ cycle.
$$(g_i\Phi=\tau (Xg_i) , \tau (Xg_i R),  \tau (Xg_i R^2), \cdots ,\tau(X g_iR^{r-1})).\eqno(7.8)$$

  \smallskip Note that $g_i (e^{\pi i/q})$ is the odd vertex of the special
   triangle $g_i\Phi$ and that $R(e^{\pi i/q})  = e^{\pi i/q}$.
   For each $Xg_i R^j$, by (7.2), there exists a unique pair $(g_j\Phi, x_{ij})
   \in \Omega _X\times X$ such that
   $x_{ij}g_iR^j\Phi = g_j\Phi\in \Omega_X$.
    By (7.1), $ \tau (Xg_i R^j )= x_{ij}g_iR^j\Phi$.
     The odd vertex of  $x_{ij}g_iR^j\Phi $ is
       $x_{ij} g_iR^j( e^{\pi i/q})=  x_{ij} g_i (e^{\pi i/q})      \in M_X$.
        Since $x_{ij}\in X$, $g_i (e^{\pi i/q})\in M_X$ and every
         odd vertex has a unique $X$-image in $M_X$, we conclude that
          $x_{ij} g_i (e^{\pi i/q})=g_i (e^{\pi i/q})$.
           Hence the odd vertex of  $ \tau (Xg_i R^j )= x_{ij}g_iR^j\Phi$
            is the same as the odd vertex of $g_i\Phi$ for every $j$.
        Hence  the special triangles in (7.8) are special
          triangles of
          a tile $\Delta$  whose odd vertex is $g_i(e^{\pi i/q})$
           ($\Delta$ is either a $q$-gon or an $r$-cluster).
           Since $R$ acts as a counter-clockwise
           rotation  about $e^{\pi i/q}$,
            the members in (7.8) are ordered following the  orientation of $M_X$.
          The following is clear.

          \smallskip
          \noindent {\bf Lemma 7.5.} {\em
         Let $\Delta$ be a tile of $M_X$ $(\Delta$ is either a $q$-gon or an $r$-cluster$)$.
         Denoted by $c_{\Delta}$
          the cycle of special triangles of $\Delta$  ordered
           in the counter-clockwise manner. The permutation representation of $R$ on $\Omega_X$ is
            $$f(R) = \prod c_{\Delta}.\eqno(7.9)$$}

     \noindent {\bf Example 7.6.} Let $M_X$ be given as in (7.5).
         Then
     $f(R) = (1, 6,\overline 2)(2, 5, \overline 3)(3, \overline 5, 4)
     (\overline 1)(4).$

  \section{ Normaliser of $X$ in $PSL(2, \Bbb  R)$}
  Let $G_X = \left <f(S), f(R)\right >$. Then $G_X$ is a subgroup of $S_n$,
     where $[G_q : X] = n$.
      The main purpose of this section is to
       determine the normaliser $N(X)$  of $X$ in $PSL(2, \Bbb R)$.
       In particular, we are able to determine whether $X$ is a normal subgroup of
         $G_q$.

\smallskip
 \noindent
 {\bf Proposition 8.1.} {\em Let $X$ be a subgroup of $G_q$ of index $n$.
   Then $N_{G_q}(X)/X
   \cong C_{S_n}(G_X)$. $X$ is a normal subgroup of $G_q$
    if and only if the order of $G_X$ is $ [G_q : X] = n$.}

 \smallskip
 \noindent {\em Proof.}
 Consider the action of $G_X$ on $\Omega _X$. It is clear that the action is transitive.
  Let $G_0$ be a one point stabiliser. By Lemma A1 of Appendix A,
  $C_{S_n}(G_X) \cong N_{G_X}(G_0)/G_0$.
  By Lemma 7.1, the action of $G_q$ on $G_q/X$ is
  isomorphic
   to the action of $G_X$ on $\Omega _X$.
Hence  $ N_{G_q}(X)/X\cong N_{G_X}(G_0)/G_0 \cong C_{S_n}(G_X)  $.

  \smallskip

   Consider the action of $G_q$ on $G_q/X$, one has
    $X \triangleleft \,G_q$ if and only if $Y = \cap \,gXg^{-1} =X$ if and only
    if $ G_X \cong G_q/Y \cong G_q/X$
    if and only if $G_X$ has order $n$.\qed

\smallskip
\noindent {\bf Discussion 8.2.} (i)
 Since $f(S)$ and $f(R)$ can be obtained easily by studying the special
  polygon $M_X$ and the order of $G_X$ can be determined by {\bf GAP}, whether
   $X$ is normal becomes an easy issue.
 In the case $q\ne 3,4,6$,  Margulis' characterisation of arithmeticity in terms
      of the commensurator implies that
 $N(X) = N_{G_q}(X)$. Hence
  $C_{S_n}(G_X) \cong N(X)/X$.
(ii) By Proposition 8.1, the group in (i) of subsection 5.1 is normal.

\smallskip

 \noindent {\bf Example 8.3.}
 Let
$$M_X = \{  {-\infty }_{_{_{\smile}} }\  \hspace{-.37cm}_{ _{_{_{_{_{_{_{1}}}}}}} }  {0}
_{_{_{\smile}} }\  \hspace{-.37cm}_{ _{_{_{_{_{_{_{2}}}}}}}}
1
 _{_{_{\smile}}}\  \hspace{-.36cm} _{_{_{_{_{_{_{_{2}}}}}}}}
 3/2
 _{_{_{\smile}}}\  \hspace{-.36cm} _{_{_{_{_{_{_{_{3}}}}}}}}
  2
    _{_{_{\smile}}}\  \hspace{-.36cm} _{_{_{_{_{_{_{_{3}}}}}}}}
 {3} _{_{_{\smile}}} \  \hspace{-.36cm} _{_{_{_{_{_{_{_{1}}}}}}}}
 {\infty}
  \}
   \eqno(8.1)$$
be the Hecke-Farey symbol of $X \subseteq G_3 = PSL(2, \Bbb Z) $ of index 12
(see Figure 5). $M_X$ has 6 even lines. They are
 $$ L_1 = (0, \infty),L_2 = (0, 1),
L_3 = (3/2, 2), L_4= (1,\infty),  L_5=(2, \infty),
 L_6=(1,2).\eqno(8.2)$$
Note that $(3,\infty)$ is paired  with $L_1$,
 $(1, 3/2)$ is paired with $L_2$, and that
$(2,3)$ is paired with $L_3$ (see (8.1) for the side pairings).
 One sees easily that all the even lines are shared by 2 special triangles.
   Following Discussion  7.3, the special triangles are labeled as in Figure 5.
 By Lemmas 7.2 and 7.5, we have
  $$f(R)  =
 (1,2, \bar 4)(  4,6, \bar 5)( 5, \bar 3,\bar 1 )(\bar 6,\bar 2, 3),\,
 f(S)= (1,\bar 1) (2,\bar 2)
(3,\bar 3) (4,\bar 4)
(5,\bar 5) (6,\bar 6).\eqno(8.3)$$

\noindent One has $|G_X|= |\left < f(S), f(R)\right >|
 = [G_3 : X] = 12. $ By Proposition 8.1, $X$ is a  normal.

\allowdisplaybreaks

\begin{center}
\begin{figure}
\beginpicture

\setcoordinatesystem units <5pt,5pt>

\setplotarea x from  -17 to 63, y from -1 to 18

\circulararc -180 degrees from 0 0 center at 7.5 0
\circulararc -180 degrees from 15 0 center at 22.5 0
\circulararc -180 degrees from 15 0 center at 18.75 0
\circulararc -180 degrees from 22.5 0 center at 26.25 0
\circulararc -180 degrees from 30 0 center at 37.5 0

\setlinear \plot 0 0 45 0 /

\setlinear \plot  0   0 0 16 /
\setlinear \plot 15  0 15 16 /
\setlinear \plot 30  0 30 16 /
\setlinear \plot 45  0 45 16 /

\put {\tiny $1$} at             2 12
\put {\tiny $\overline 1$} at  43  12
\put {\tiny $4$} at            17 12
\put {\tiny $\overline 4$} at   13 12

\put {\tiny $5$} at           32  12
\put {\tiny $\overline 5$} at  28  12

\put {\tiny $2$} at  7.3 9
\put {\tiny $\overline 2$} at 18.75 5
\put {\tiny $3$} at  26.25 5
\put {\tiny $\overline 3$} at  37.5  9
\put {\tiny   $6$} at  22.5  9
\put {\tiny $\overline 6$} at  22.5  6



\put {$\scriptsize{0}$} at 0 -3
\put {$\scriptsize{1}$} at 15 -3
\put {$\scriptsize{2}$} at 30 -3
\put {$\scriptsize{3/2}$} at 22.5 -3
\put {$\scriptsize{3}$} at 45 -3

\put {\bf Figure 5} at 22.5 -7
\endpicture
\end{figure}
\end{center}

\vspace{-.4cm}

\smallskip

 \noindent {\bf Example 8.4.} Let $\lambda = \lambda_4 = \sqrt2$ and let
$$M_X = \{  {-\infty }_{_{_{\smile}} }\  \hspace{-.37cm}_{ _{_{_{_{_{_{_{1}}}}}}} }
 {0}
_{_{_{\smile}} }\  \hspace{-.37cm}_{ _{_{_{_{_{_{_{2}}}}}}}}
1/\lambda
 _{_{_{\smile}}}\  \hspace{-.36cm} _{_{_{_{_{_{_{_{1}}}}}}}}
 \lambda
 _{_{_{\smile}}}\  \hspace{-.36cm} _{_{_{_{_{_{_{_{3}}}}}}}}
  3/\lambda \,
    _{_{_{\smile}}}\  \hspace{-.36cm} _{_{_{_{_{_{_{_{2}}}}}}}}
 {2\lambda}\,
  _{_{_{\smile}}} \  \hspace{-.36cm} _{_{_{_{_{_{_{_{3}}}}}}}}
 {\infty}
  \}
  \eqno(8.4)$$
be the Hecke-Farey symbol of $X\subseteq G_4$  of index 8 (see Figure 6).
$M_X$ has four even lines, $L_1=(0,\infty)$, $L_2= (0, 1/\lambda)$,
  $L_3= (\lambda, 3/\lambda)$, and $L_4= (\lambda, \infty)$.
One sees  from Figure 6 that $$f(R)  =
 (1,2, \bar 1, \bar 4)(  4,3, \bar 2, \bar3),\,
f(S) =
(1,\bar 1) (2,\bar 2)
(3,\bar 3) (4,\bar 4).\eqno(8.5)$$

 \noindent   One has  $|G_X| = 16 >  8  = [G_4 : X]$.
   By Proposition 8.1,  $X$ is not normal.
\smallskip

\noindent {\bf Example 8.5.} Let $P = \{x_1=0, x_2, \cdots, x_{q}= \infty\}$
 be the depth 1 $q$-gon in the right half plane and let  $r$ $(1\le r <q$) be a divisor of $q$. Let
   $$ M_A = \{
  -\infty \,\, _{_{_{\smile}}} \hspace{-.22cm}_{ _{_{_{_{_{_{_{\circ}}}}}}} }
   x_1\,\,
  _{_{_{\smile}}} \hspace{-.22cm}_{ _{_{_{_{_{_{_{\circ}}}}}}} }
  x_2 \,\, _{_{_{\smile}}} \hspace{-.22cm}_{ _{_{_{_{_{_{_{\circ}}}}}}} }
  \cdots
 x_{q-1}  \,\,\,\, _{_{_{\smile}}} \hspace{-.22cm}_{ _{_{_{_{_{_{_{\circ}}}}}}} }
   x_{q}\} \mbox{ and }
  M_B = \{
  -\infty \,\, _{_{_{\smile}}} \hspace{-.22cm}_{ _{_{_{_{_{_{_{\circ}}}}}}} }
   x_1\,\,
  _{_{_{\smile}}} \hspace{-.22cm}_{ _{_{_{_{_{_{_{\circ}}}}}}} }
  x_2 \,\, _{_{_{\smile}}} \hspace{-.22cm}_{ _{_{_{_{_{_{_{\circ}}}}}}} }
  \cdots
 x_r  \,\,\,\, _{_{_{\smile}}} \hspace{-.37cm}_{ _{_{_{_{_{_{_{e_r(1)}}}}}}} }
   x_{q}\}
    .$$
  $f(S)$ is the identity permutation for both $A$ and $B$.
   $f(R)$ is a $q$-cycle (resp. $r$-cycle) for $A$ (resp. $B$).  By Proposition
     8.1, both $A$ and $B$ are normal subgroups of $G_q$.

\section{Hurwitz-Nielsen realisation problem for normal subgroups}
 A transitive subgroup $G$ of $S_d$ is called {\em regular} if every non-identity
  element of $G$ is fixed point free.
The result in this section
 reveals the connection between normal subgroups of index $d$
  of $G_q$ and regular subgroups of $S_d$.

\subsection
{Geometric invariants revisited }
Let $M_X=(P, I_X)$ be a special polygon of $X$.
The geometric invariants of $X$ in terms of $f(S)$ and $f(R)$ are given as follows.

\begin{enumerate}
\item[(i)] $I_X$ has
 $ \tau_2$ elliptic elements of order 2 that are conjugates of $S$, where $\tau_2 $ is the number  of one cycles of $f(S)$.
 \item [(ii)] Let $r$ ($1 < r\le q$) be a divisor of $q$.
 $I_X$ has
 $ v_r$ elliptic elements of order $r$ that are conjugates of $R^{q/r}$, where $v_r $ is the number  of $q/r$-cycles of $f(R)$.
 \item[(iii)] $P$ has $v_{r_i} $ $q/r_i$-clusters.
 The number $n_0$ of  $q$-gons of $M_X = (P, I_X)$
 is the number of $q$-cycles of $f(R)$.
  The index $[G_q :X]$ is given by
   $n_0q + \sum v_r q/r$.

 \item[(iv)] The number of cusps is the number of cycles of $f(T)$,
  where $f(T)$ is the permutation representation of $T$ on $\Omega_X$
  (see (v) of
  subsection 6.1).

\item[(v)] $I_X$ has
  $ f = (n_0(q-2) + \sum_{r|q} v_r(q/r-2) + 2- \tau_2)/2$
   side pairings of infinite order.

\end{enumerate}

\smallskip
\noindent {\em Proof.} By Lemma 7.5,
it is clear  that (a) an $s$-cycle $(s <q)$ of $f(R)$ corresponds
 to an $s$-cluster of $P$ as well as a side pairing which is a conjugate of $R^{s}$,
  (b) a $q$-cycle of $f(R)$ corresponds to a $q$-gon of $P$.
  Similarly, each one cycle of $f(S)$ corresponds to a side pairing which is
   a conjugate of $S$.
(i)-(iv) of the above is clear. See Lemma 6.1 for (v). \qed

\allowdisplaybreaks

\begin{center}
\begin{figure}
\beginpicture

\setcoordinatesystem units <6pt,6pt>

\setplotarea x from  -18.5 to 63, y from -1 to 14.5

\circulararc -180 degrees from 0 0 center at 3.75 0
\circulararc -180 degrees from 7.5 0 center at 11.25 0

\circulararc -180 degrees from 15 0 center at 18.75 0
\circulararc -180 degrees from 22.5 0 center at 26.25 0


\setlinear \plot 0 0 30 0 /

\setlinear \plot  0   0    0 13 /
\setlinear \plot 15  0    15 13 /
\setlinear \plot 30  0    30 13 /

\put {$1$} at             2 9
\put {$\overline 4$} at   13 9
\put {$4$} at             17 9
\put {$\overline 3$} at    28  9

\put {$ 3$} at 18.75 5
\put {$\overline  2$} at  26.25 5

\put {$\overline 1$} at 11.25 5
\put {$ 2$} at  3.75 5


\put {$\scriptsize{0}$} at 0 -3
\put {$\scriptsize{1/\lambda}$} at 7.5 -3

\put {$\scriptsize{\lambda }$} at 15 -3
\put {$\scriptsize{2\lambda}$} at 30 -3
\put {$\scriptsize{3/\lambda}$} at 22.5 -3


\put {\bf Figure 6} at 16 -7
\endpicture
\end{figure}
\end{center}

\vspace{-.7cm}

\subsection{Geometric invariants for normal subgroups}
 Suppose that  $X$ is normal of index $d$.
 By Proposition 8.1,
  $G_X = \left < f(R), f(S)\right > $
   is a transitive subgroup of order $d$ of the
    symmetric group $S_d$.
   Since   every nonidentity
   element acts freely,
     $f(R)$ is a product of
     $d/u$ $u$-cycles for some $u|q$. Similarly, $f(S)$ is either 1 or a product
      of $d/2$ 2-cycles.
        The geometric invariants can be determined by subsection 9.1. In particular,

\begin{enumerate}

\item[(vi)] $I_X$ has
  $  f = d(u-2)/u + 2- \tau_2)/2$
   side pairings of infinite order $(\tau_2=0$ if $S\notin X$).

\end{enumerate}

  The following fact is well known. Unlike most of the results in the
   literature,  our study
  gives the set of generators of such $X$'s in matrix forms (see Example 8.5).

 \smallskip

 \smallskip

 \noindent {\bf Lemma 9.1.} {\em
  Let $A$ and $B$ be given as in Example $8.5$
   and
 let $X$ be a normal subgroup of index $d$ of
  $G_q$ that contains  $S$. Then $X$ is isomorphic to either $A$ or $B$.}

\smallskip
\noindent {\em Proof.}
Let $M_X= (P, I_X)$ and let $G_X=\left < f(S), f(R)\right >$.
 Since $X$ is normal, $G_X$ is regular of order $d$.
 Since $S\in X$ and $X$ is normal, $f(S) = 1$ (see (i) of
  subsection 9.1).
  Hence $G_X = \left <f(R)\right >$.
   Since $G_X$ is transitive,  $f(R)$ must be  a $d$-cycle.
   It follows that  $P$ is either a $q$-gon (if $d=q$) or a $d$-cluster $\Phi_d$ (if $d
    < q)$.
    As $S\in X$ and $X$ is normal,  the even lines of $P$
     are self paired by conjugates of $S$. This completes the proof of the proposition.\qed

\subsection{Realisation of groups as normal subgroups  }
 Throughout the subsection,
  $Y$ is a free product of a free group $F_f$ of rank $f$, $\pi_2$
    copies of $\Bbb Z_2$, and $v_{r_i} $ copies of $\Bbb Z_{r_i}$, where
    $r_i \in \Delta_0 =\{ r\,:\, r|q, \, 3\le r\le q\}.$
Proposition  9.2 is an immediate consequence of Lemma  9.1.

     \smallskip
     \noindent {\bf Proposition 9.2.} {\em Let $q\ge 3$.
      Suppose that $\pi_2 \ne0$. Then
       $Y$ can be realised as a normal subgroup of $G_q$
       that contains $S$
       if and only if $Y$ is a free product of $q$ copies of $\Bbb Z_2$
        or a free product of $\Bbb Z_{q/r}$ and $r$ copies of $\Bbb Z_2$,
         where $r|q$ and $r < q$.

       }


  \smallskip
     \noindent {\bf Proposition 9.3.} {\em Let $q\ge 3$.
      Suppose that $\pi_2 \ne0$. Then $Y$ can be realised as a normal subgroup
     of index $d$  of $G_q$
       that does not contain $S$
       if and only if }(i) {\em
 $Y$ is a free product of  $2d/q $ copies of $\Bbb Z_{2}$
 and a free group of rank $f$, where $2f= 2d(q/2-2)/q +2$, } (ii) {\em
 $S_d$ has a regular subgroup $G = \left < \alpha, \beta\right > $ of order $d$,
    where $\alpha$ is a product of $2$-cycles and $\beta$ is a product
     of $q/2$-cycles.
     In particular, $q$ is even.
     Note that $d$ is determined by $f$ and $q$. }

\smallskip
\noindent {\em Proof.}
Let $Y$  and  $G$ be given  as in (i) and (ii).
 Since $G_q$ is a free product of $S$ and $R$, the map  $\chi \,:\, G_q \to
 G $
 defined by $\chi(S) = \alpha, \chi(R) = \beta$ is a homomorphism.
  Let $X$ be the kernel of $\chi$ and let $M_X = (P, I_X)$
   be a special polygon of $X$. One sees   that  $f(S)$ (resp.
    $f(R)$) and $\alpha$ (resp. $\beta$) have the same cycle decomposition.
   Hence $f(S)$ is a product of $2$-cycles,
 $f(R)$ is a product of   $q/2$-cycles, and
  $P$ is a union of  $q/2$-gons ($2d/q$ of them). By (vi) of subsection 9.2,
  $2f = 2d(q/2-2) /q +2$. Hence
 $I_X$ has $2d/q$ elements that are conjugates of $R^{q/2}$
       and $f$ elements of infinite order. In particular, $Y \cong X \triangleleft G_q$.

  Conversely, suppose that $Y\cong X\triangleleft G_q$, where $[G_q : X] = d$,
   $S\notin X$. Let $G_X=\left <f(S), f(R)\right >$.
   By (i)-(vi) of subsections 9.1 and 9.2, $G_X$ is regular, $f(S)$
    (resp $f(R)$) is a product of $2$-cycles (resp. $q/2$-cycles).  Further,
   $X$ s a free product
    of $2d/q$ copies of $\Bbb Z_2$ and a free group $F_f$, where $2f= 2d(q/2-2)/q+2$.
   Hence (i) and (ii) holds.\qed

\smallskip
The following proposition can be proved by applying the proof of Proposition 9.3.

\smallskip
\noindent {\bf Proposition 9.4.} {\em Let $q\ge 3$. Suppose that $\pi_2 =0$.
  Then $Y$ can be realised as
    a normal subgroup of index $d$  of $G_q$ if and only if
 \begin{enumerate}
   \item[(i)] $Y$ is a free group of rank $f$, where $2f = d(q-2)/q +2$,
   $S_d$ has regular subgroup $\left < \alpha, \beta\right > $ of order $d$,
    where $\alpha$ is a product of $2$-cycles and $\beta$ is a product
     of $q$-cycles,  or
   \item[(ii)]
 $Y$ is a free product of  $d/r $ copies of $\Bbb Z_{q/r}$
 and a free group of rank $f$, where $2f= d(r-2)/r +2$,
 $S_d$ has regular subgroup $\left < \alpha, \beta\right > $ of order $d$,
    where $\alpha$ is a product of $2$-cycles and $\beta$ is a product
     of $r$-cycles. Note that $d$ is  determined by $f, q$ and $r$.
 \end{enumerate}}

  \section {Construction of all maps on compact orientable surfaces }
\subsection{Known results}
A {\em map}  $\mathfrak M $ on a compact orientable surface $\mathcal X$ is an embedding of
 a finite connected graph $\mathcal G$ in $\mathcal  X$ such that
  the connected components (faces of $\mathfrak M$) of $\mathcal X \setminus \mathcal G$
are simply connected.
An edge $e$ of $\mathfrak M$ is called a {\em segment} if $e$ has two vertices.
 An edge  homeomorphic to the circle $S^1$ is called a {\em loop}.
 An edge homeomorphic to $[0,1]$ with only one vertex is called a {\em free edge}.
  Segments and loops are also called {\em non-free edges}.
  The {\em darts} of $\mathfrak M$ are {\em directed edges}.
    Each non-free edge gives two darts whereas a free edge
    gives only one dart (see pp. 276 of [JS]).
    The two darts associated with a non-free edge $e$ travel along $e$ in opposite
     directions and the only dart associated with a free edge $e$ (with vertex $v$)
      is the directed edge that points towards $v$.
     Denoted by $\Omega$ the set of darts
  of $\mathfrak M$.
$\mathfrak M$  can be characterised completely by two
 permutations $r_1$ and $r_2$ on $\Omega$. $r_1$ is the permutation that (i) fixes
  the darts associated with the free edges,  and  (ii) transposes  the two darts
  associated with the  non-free edges.
  $r_2$ is the  permutation whose cycles correspond to the faces of $\mathfrak M$ (following the orientation of $\mathcal X$).
 Consequently,  the map $\mathfrak M$ can be represented by
  $$(G, \Omega, r_1, r_2),\eqno(10.1)$$
  where $G = \left < r_1, r_2\right >$.
    Suppose that the order of $r_2$ is $q$. Then
    $\theta \,:\, G_q \to G= \left <r_1, r_2\right > \subseteq S_{\Omega}$
     defined by $\theta (S)=r_1$ and  $\theta (R) = r_2$ is a homomorphism.
      Since the graph $\mathcal G$ is connected, $\theta$ induces a transitive action of
       $G_q$ on $\Omega$ via $A(w) = \theta (A) (w)$, for all $A\in G_q$, $w\in \Omega$.
        Let $X$ be a one point stabiliser of the action. It is clear that
         the action of $G_q$ on $G_q/X$ is isomorphic to the action of $G$ on $\Omega$
          and that $(G, \Omega, r_1, r_2)$
           is isomorphic to
           $(G_q/X_0, G_q/X, S, R)$, where $X_0= \cap _{g\in G_q} gXg^{-1}$
            (see pp. 283 of [JS]). Equivalently,
     $$\mathfrak M\cong (G, \Omega, r_1, r_2)\cong (G_q/X_0, G_q/X, S, R).\eqno(10.2)$$        $X$ is known as a {\em map subgroup} of $\mathfrak M$.
           Since the action is transitive, the map subgroups are conjugate to each other.
         See pp.63 of  [CS] for more detail.

\subsection{Maps associated with subgroups of $G_q$}  Let $M_X = (P, I_X) = X\setminus \Bbb H$ be a special polygon  of $X\subseteq G_q$ given as in Section 4. Define $M(X)$ as follows.

\begin{enumerate}
\item[(i)] $V = $ set of vertices of $M(X) =$ the set of equivalence classes of cusps of $M_X$,
\item[(ii)] $E=$ set of edges of $M(X) =$   the
  set of  equivalence classes of even lines of $M_X$,
  \item[(iii)] $F=$ set of faces of $M(X)=$ the set of
    clusters and $q$-gons of  $M_X$.

    \end{enumerate}

   As the faces of $M(X)$ are simply connected,  $M(X) =(V, E, F)$
  is a  map.

 \smallskip
 \noindent {\bf Lemma 10.1.} {\em
  Suppose that $[G_q : X]< \infty$. Then $ M(X) \cong (G_X, \Omega_X, f(S), f(R)).$}

 \smallskip
 \noindent {\em Proof.} Let
 $E_1$ and $E_2$ be given as in Discussion  7.3.  Let $e\in E_1$. Then $e$ is not paired with itself by
    elements of order 2. Hence $e$ is either a loop or a segment.
   Consequently,
    $e$ is a non-free edge.
    In the case $e\in E_2$, $e$ is paired with itself by an element of order 2. Hence
     $e$ is homeomorphic to $[0,1]$ with only one vertex.
     It follows that $e$ is a free edge.
    Note that every non-free edge $e$ is shared by two special triangles  and
     every free edge $e$ belongs to exactly one special triangle.
          This allows us to construct a one to one
      correspondence $\nu$ between the set of special triangles $\Omega_X$
       and the set of darts $\Omega$ such that
        $\nu(f(S)(g_i\Phi)) = r_1(\nu(g_i\Phi))$ and that
    $\nu(f(R)(g_i\Phi)) = r_2(\nu(g_i\Phi))$. Hence the following is a
     commutative diagram.
    $$
    \begin{array}{ccl }
       \Omega_X \times  G_X  &  \longrightarrow &  \Omega_X\\
(\nu, \chi  ) \phantom {\Bigg |}\downarrow\,\,\,\,\,\,\,\,\, && \, \downarrow  \nu\\
  \, \Omega \times  G  & \longrightarrow  &  \,\Omega \\
     \end{array} \eqno(10.3)$$
 where the horizontal arrows represent the group actions and    $\chi (f(S))) = r_1$,
  $\chi(f(R)) = r_2$. As a consequence,
       $ M(X)
   \cong (G, \Omega, r_1,r_2)
  \cong (G_X, \Omega_X, f(S), f(R)).$\qed

\smallskip
\noindent {\bf  Discussion 10.2. }
 Let $\mathfrak M$ be a map.
 One knows very little about the   map subgroups of $\mathfrak M$ except that they
 are the one point stabilisers. Conversely, for each subgroup $X$ of $G_q$,
  one cannot really {\em visualise} $(G_q/X_0, G_q/X, S, R)$ as a map as it is
   not very easy to describe the incidence relations  of this map.
  The following proposition implies that
    every map (with $o(r_2)=q$) takes the
    form $M(X)$, where $X$ is a subgroup of finite index of $G_q$.
     Both the map $M(X) $ and its map subgroups can be described
    explicitly as
      $M(X)$ can be  described by its special polygon and the map
       subgroups of $M(X)$, which are conjugates of $X$, can be
        described by the set of independent generators $I_X$.

 \smallskip
 \noindent {\bf Proposition 10.3.} {\em
   Let $\mathfrak M$ be a map, where $o(r_2)=q$.
   Then $\mathfrak M
    \cong M(X)$ for some $X \subseteq G_q$. Further, the map subgroups of $M(X)$ are
     conjugates of $X$.}

 \smallskip
 \noindent {\em Proof.} By the results in subsection 10.1, $\mathfrak M$
  is isomorphic to
  $(G_q/X_0, G_q/X, S, R)$. By Lemmas 7.1, 7.2, and 7.5,
  $(G_q/X_0, G_q/X, S, R)$ is isomorphic to
   $(G_X, \Omega_X, f(S), f(R))$.
    By  Lemma 10.1, one has  $\mathfrak M\cong (G_X, \Omega_X, f(S), f(R))\cong M(X)$.

    $g\in G_q$ fixes
      $g_i\Phi\in \Omega_X$ if and only if there exists some $x\in X$
       such that $xg_i g\Phi = g_i\Phi$ (see (7.2)).
      Since the action of $G_q$ on $\Omega_X$
       is fixed point free,  $g$ fixes $g_i\Phi$ if and only if
        $ g_i^{-1} xg_i g =1$.
     This completes the proof of the proposition.\qed


 \smallskip
 \noindent {\bf Proposition 10.4.} {\em
   The automorphism group of $M(X)$ is
 $C_{S_{n}}(G_X)$, where $n=[G_q:X]$.}

  \smallskip
 \noindent {\em Proof.} Let $\sigma\in S_d$.
  The  incidence relations of $M(X)$ is  determined by $G_X
   =\left <f(S), f(R)\right >$. Hence
    $\sigma\in $ Aut$\, M(X)$ if and only if $[\sigma,G_X]=1$.
    This completes the proof of the lemma.
 \qed

\section{ Congruence Subgroup Problem for $\Gamma = PSL(2, \Bbb Z)$}
Let $\Gamma = G_3 = PSL(2, \Bbb Z)$. The {\em principal congruence subgroup} of $\Gamma$ of {\em level} $r\in \Bbb N$ is
  $$\Gamma(r) = \{ x\in \Gamma \,:\, x \equiv \pm I\,\,(mod\,\,  r)\}.\eqno(11.1)$$

\smallskip
\noindent  $X\subseteq \Gamma$ is a {\em congruence subgroup} if $\Gamma(r) \subseteq X$ for some $r$.
 Consider the action of $T$ on $\Gamma/X$,
 the order of $T $ on $\Gamma/X$  is called the {\em level} of $X$.
  By a result of Wohlfahrt [W], $X \subseteq \Gamma$ is a congruence subgroup if and only
   if $\Gamma (n)\subseteq X$.
   To the best of our knowledge, the congruence test
    developed by
         Tim Hsu [H] is the most effective one in the literature
          (see [LLT2] for another test). His test
          can be implemented as long as  the permutation
          representations of $T$ and $U$ on $\Gamma/X$ can be determined
          (Theorem 3.1 of [H]),
            where
           $$T =\left (
\begin{array}{rr}
1 & 1\\
0 &1 \\
\end{array}
\right ) ,\,
 U =\left (
\begin{array}{rr}
1 & 0\\
1 &1\\
\end{array}
\right ).\eqno(11.2)$$

Incidently, an algorithm that determines the
 permutation representations of $T$ and $U$ on $\Gamma/X$  is not included in [H].
By Lemma 7.1,  such action is isomorphic to the action of
 $T$ and $U$ on $\Omega_X$.
Following Lemmas 7.2 and 7.5, $f(T)= f(R^{-1}S)$ and $f(U) = f(RS)$ can be determined
 easily as long as $X$ is given in terms of a special polygon and
  a set of independent generators.
    As a consequence,  Hsu's algorithm can be
    implemented with ease.

\smallskip
\noindent {\bf Example 11.1.}
 Let
$$ M_X = \{  {-\infty }_{_{_{\smile}} }\  \hspace{-.37cm}_{ _{_{_{_{_{_{_{\bullet}}}}}}} }
 0
_{_{_{\smile}} }\  \hspace{-.37cm}_{ _{_{_{_{_{_{_{\circ}}}}}}}}
1
 _{_{_{\smile}}}\  \hspace{-.36cm} _{_{_{_{_{_{_{_{\circ}}}}}}}}
 2
_{_{_{\smile}} }\  \hspace{-.37cm}_{ _{_{_{_{_{_{_{\circ}}}}}}}}
3
 _{_{_{\smile}}}\  \hspace{-.36cm} _{_{_{_{_{_{_{_{\bullet}}}}}}}}
 \infty
  \}\eqno(11.3)$$

\smallskip
\noindent be the Hecke-Farey symbol of a subgroup $X$ of index 11 of $PSL(2, \Bbb Z)$
 (see Figure 7).
 Then
$ f(S)= (1, \bar 1)(2,\bar 2)(3, \bar 3)(4, \bar 4)$ and
 $ f(R)  = (\bar 1) (1,5,\bar 2)(2,6,\bar 3)(3,7,\bar 4)( 4).$
  Hence
   $$ f( T)= (\bar 1, \bar 2, \bar 3, \bar 4, 4,7,3,6,2,5,1),\,\,
  f(U) = (\bar 1,5,\bar 2, 6, \bar 3, 7, \bar 4, 4,3,2,1)
  .\eqno(11.4)$$

\smallskip
\noindent By the algorithm given in Section 3 of [H], $X$ is non-congruence.
 Note that   $\cap \,gXg^{-1}  \cong\left
  < f(T),f( U)\right > \cong A_{11}$ is the  alternating group on 11 letters.
   The  group $\cap \,gXg^{-1} $  was first studied by Magnus [M]
   as part of his study of non-congruence subgroups of $\Gamma$.

\allowdisplaybreaks

\begin{center}
\begin{figure}
\beginpicture

\setcoordinatesystem units <4pt,4pt>

\setplotarea x from   -26.5 to 83, y from -1 to 22

\circulararc -84 degrees from 0 0 center at 7.5 0

\circulararc  -90 degrees from 7.5 7.5 center at 7.5 0

\circulararc -87 degrees from 15 0 center at 22.5 0
\circulararc  -88 degrees from 22.9 7.5 center at 22.5 0

\circulararc -87 degrees from 30 0 center at 37.5 0
\circulararc -87 degrees from 38 7.5 center at 37.5 0

\circulararc -60 degrees from 45 0 center at  60 0
\circulararc  60 degrees from 0 0 center at  -15 0

\setlinear \plot 0 0 45 0 /

\setlinear \plot  0   0 0 20 /
\setlinear \plot 15  0 15 20 /
\setlinear \plot 30  0 30 20 /
\setlinear \plot 45  0 45 20 /

\setlinear \plot -7.5  13 -7.5 20 /
\setlinear \plot 52.5  13  52.5 20 /

\put {$\bullet$} at   -7.5 13
\put {$\bullet$} at  52.5  13


\put {$    2$} at   17 13
\put {$  3$} at  32  13
\put {$ 1$} at   2  13
\put {$  \overline  4$} at  43  13

\put {$ \overline 2$} at   13 13
\put {$ \overline 3$} at  28  13
\put {$ \overline 1$} at   -2  13
\put {$  4$} at  47  13

\put {$5$} at  7.3 5.5

\put {$7$} at  37.5  5.5
\put {$6$} at  22.5  5.5


\put {${ \circ } $ } at  7.5 7.5
\put {$ \circ$} at  22.5 7.5
\put {$ \circ$} at  37.5 7.5

\put {$\scriptsize{0}$} at 0 -3
\put {$\scriptsize{1}$} at 15 -3
\put {$\scriptsize{2}$} at 30 -3

\put {$\scriptsize{3}$} at 45 -3

\put {\bf Figure 7} at 22.5 -7
\endpicture
\end{figure}
\end{center}

\vspace{-.5cm}

\section* {Appendix A  }

In this appendix,  $\Omega$ is a finite set, $S_{\Omega}$ is the
 symmetric group on $\Omega$, and $G \subseteq S_{\Omega}$.




\smallskip
\noindent {\bf Lemma A1.} {\em
Suppose that $G$ acts  transitively on $\Omega$.
      Then $C_{S_{\Omega}}(G) \cong N_G(G_d)/G_d$,
       where $G_d$ is the one point stabiliser of $d \in \Omega$.}

       \smallskip
       \noindent {\em Proof.} Since $G$ is transitive,
       The action of $G $ on $\Omega$ is isomorphic to the   action of $G$ on the set of
        cosets $G/G_d$. Without loss of generality, we may assume that
        $\Omega = G/G_d$ and that $x(G_d) = xG_d$ for $x \in G$.
            Let $x\in C_{S_{\Omega}}(G)$. Then $x(G_d) = e_xG_d\in G/G_d$ for some $e_x \in G$.
     For each $g \in G_d$, one has $gx(G_d) = xg(G_d) = x(G_d)$.
     Hence $ge_xG_d = e_x G_d$. This implies that
     $e_xG_d \in N_G(G_d)/G_d$. As a consequence, one can show that
     $C_{S_{\Omega}}(G) \cong N_G(G_d)/G_d$ by studying the
      homomorphism  $\Phi\, :\, C_{S_{\Omega}}(G) \to N_G(G_d)/G_d$ defined
       by $\Phi(x) = e_x^{-1} G_d$. Note that for each $  r \in N_G(G_d)$,
        the permutation defined by $\sigma (g G_d) = gr^{-1} G_d$ commutes with $G$
         which implies that $\Phi$ is surjective  ($e_{\sigma} = r^{-1}$ and
          $\Phi (\sigma) = r^{}G_d$).
    \qed




\bigskip{\small

\noindent Cheng Lien Lang\\
\noindent Department of Mathematics, I-Shou  University, Kaohsiung, Taiwan.

\noindent   \texttt{cllang@isu.edu.tw}


\smallskip
\noindent Mong Lung Lang \\
\noindent Singapore 669608,
Singapore.

\noindent \texttt{lang2to46@gmail.com}}

\medskip


\medskip

\end{document}